\documentclass[12pt,fleqn,leqno]{article}
\usepackage{amsmath,amssymb}

\usepackage[dvips]{graphicx}
\usepackage[dvips]{color}

\allowdisplaybreaks

\numberwithin{equation}{section}

\makeatletter
\renewcommand{\section}{\@startsection{section}{1}{0pt}{20pt}{6pt}{\large\bf}}
\renewcommand{\@seccntformat}[1]{\csname the#1\endcsname.\ }

\def\footnoterule{\kern -3pt \hrule width 2.7 true cm \kern 2.6pt}

\def\h{\hspace}

\def\ni{\noindent}
\def\p{\!+\!}
\def\m{\!-\!}

\def\EE{\mathsf E\:\!}
\def\PP{\mathsf P}

\def\cF{{\cal F}}

\def\R{I\!\!R}
\def\LL{I\!\!L}
\def\eps{\varepsilon}

\begin{document}

\title{\bf Optimal Real-Time Detection of \\ a Drifting Brownian Coordinate}
\author{P. A. Ernst, G. Peskir \& Q. Zhou}
\date{}
\maketitle

%%%%%%%%%%%%%%%%%%%%%%%%%%%%%%%%%%%%%%%%%%%%%%%%%%%%%%%%%%%%%%%%%%%%%%%%%%%%%%%
%%% Research Report %%%
%%%%%%%%%%%%%%%%%%%%%%%%%%%%%%%%%%%%%%%%%%%%%%%%%%%%%%%%%%%%%%%%%%%%%%%%%%%%%%%

%%% {\par \vspace{-14pc} \noindent \footnotesize
%%%
%%% \ni \emph{Research Report} No.\ 1, 2018, \emph{Probab.\ Statist.\
%%% Group Manchester} (37 pp)
%%%
%%% \vspace{13pc} \par}

%%%%%%%%%%%%%%%%%%%%%%%%%%%%%%%%%%%%%%%%%%%%%%%%%%%%%%%%%%%%%%%%%%%%%%%%%%%%%%%
%%% Abstract %%%
%%%%%%%%%%%%%%%%%%%%%%%%%%%%%%%%%%%%%%%%%%%%%%%%%%%%%%%%%%%%%%%%%%%%%%%%%%%%%%%

{\par \leftskip=2.7cm \rightskip=2.7cm \footnotesize

Consider the motion of a Brownian particle in three dimensions,
whose two spatial coordinates are standard Brownian motions with
zero drift, and the remaining (unknown) spatial coordinate is a
standard Brownian motion with a non-zero drift. Given that the
position of the Brownian particle is being observed in real time,
the problem is to detect as soon as possible and with minimal
probabilities of the wrong terminal decisions, which spatial
coordinate has the non-zero drift. We solve this problem in the
Bayesian formulation, under any prior probabilities of the non-zero
drift being in any of the three spatial coordinates, when the
passage of time is penalised linearly. Finding the exact solution to
the problem in three dimensions, including a rigorous treatment of
its non-monotone optimal stopping boundaries, is the main
contribution of the present paper. To our knowledge this is the
first time that such a problem has been solved in the literature.

\par}

%%%%%%%%%%%%%%%%%%%%%%%%%%%%%%%%%%%%%%%%%%%%%%%%%%%%%%%%%%%%%%%%%%%%%%%%%%%%%%%
%%% MSC & Key words %%%
%%%%%%%%%%%%%%%%%%%%%%%%%%%%%%%%%%%%%%%%%%%%%%%%%%%%%%%%%%%%%%%%%%%%%%%%%%%%%%%



%%%%%%%%%%%%%%%%%%%%%%%%%%%%%%%%%%%%%%%%%%%%%%%%%%%%%%%%%%%%%%%%%%%%%%%%%%%%%%%
%%% Sections %%%
%%%%%%%%%%%%%%%%%%%%%%%%%%%%%%%%%%%%%%%%%%%%%%%%%%%%%%%%%%%%%%%%%%%%%%%%%%%%%%%

\vspace{12pt}

\section{Introduction}
%%%%%%%%%%%%%%%%%%%%%%

Imagine the motion of a Brownian particle in three dimensions, whose
two spatial coordinates are standard Brownian motions with zero
drift, and the remaining (unknown) spatial coordinate is a standard
Brownian motion with a non-zero drift. Given that the position $X$
of the Brownian particle is being observed in real time, the problem
is to detect as soon as possible and with minimal probabilities of
the wrong terminal decisions, which spatial coordinate has the
non-zero drift. The purpose of the present paper is to derive the
solution to this problem in the Bayesian formulation, under any
prior probabilities of the non-zero drift being in any of the three
spatial coordinates, when the passage of time is penalised linearly.

The loss to be minimised over sequential decision rules is expressed
as the linear combination of the expected running time and the
probabilities of the wrong terminal decisions. This problem
formulation of sequential testing dates back to \cite{Wa} and has
been extensively studied to date (see \cite{JP} and the references
therein). The linear combination represents the Lagrangian and once
the optimisation problem has been solved in this form it will also
lead to the solution of the constrained problem where upper bounds
are imposed on the probabilities of the wrong terminal decisions.
The central focus of the present paper is on the Lagrangian and the
methods needed to solve the problem in this form. The constrained
problem itself will not be considered in the present paper as this
extension is somewhat lengthy and more routine.

Standard arguments show that the initial optimisation problem can be
reduced to an optimal stopping problem for the posterior probability
process $\varPi$ of the non-zero drift being in the spatial
coordinates given $X$. A canonical example of $X$ in one dimension
is Brownian motion having one among two constant drifts (see
\cite{Mi} and \cite{Sh-1}). In this case $\varPi$ is a
one-dimensional Markov/diffusion process. This problem has also been
solved in finite horizon (see \cite{GP}). Books \cite[Section
4.2]{Sh-2} and \cite[Section 21]{PS} contain expositions of these
results and provide further details and references. Signal-to-noise
ratio in these problems (defined as the difference between the two
drifts divided by the diffusion coefficient) is constant. Sequential
testing problems for $X$ in one dimension where the signal-to-noise
ratio is not constant were studied more recently in \cite{GS} and
\cite{JP}. In these problems $\varPi$ is no longer Markovian,
however, the process $(\varPi,X)$ is a two-dimensional
Markov/diffusion process with the infinitesimal generator of
\emph{parabolic} type.

Another canonical example of $X$ in one dimension is Brownian motion
having one among three or more constant drifts (see \cite{SW} for a
discrete time analogue). This problem has been studied more recently
in \cite{ZS} (see also \cite{DPS} for a Poisson process analogue).
The Markov/diffusion process $\varPi$ is two-dimensional and its
infinitesimal generator is also of \emph{parabolic} type.

Related sequential testing problems for $X$ in three or more
dimensions when each coordinate process of $X$ can have a non-zero
drift have been studied in \cite{LPXG} and \cite{BK}. These problems
contain an element of optimal control as well in deciding which
coordinate process should be observed at any given time. The former
paper contains a review of other related papers (such as \cite{PR})
and the latter paper shows that the Markov/diffusion process
$\varPi$ is one-dimensional even if one admits infinitely many
coordinate processes of $X$ in the problem formulation.

In contrast to all the sequential testing problems studied to date
we will see below that the two-dimensional Markov/diffusion process
$\varPi$ in the sequential testing problem of the present paper has
the infinitesimal generator of \emph{elliptic} type. Moreover, we
will also see that the optimal stopping boundaries are
\emph{non-monotone} as functions of the coordinate variables. This
fact itself presents a formidable challenge as to our knowledge no
rigorous treatment of non-monotone optimal stopping boundaries has
been exposed in the probabilistic literature as yet. Finding the
exact solution to the problem for $X$ in three dimensions, including
a rigorous treatment of its non-monotone optimal stopping
boundaries, is the main contribution of the present paper. To our
knowledge this is the first time that such a problem has been solved
in the literature.

\section{Outline of the paper}
%%%%%%%%%%%%%%%%%%%%%%%%%%%%%%

The exposition of the material is organised as follows. In Section 3
we derive the optimal stopping problem for
$\varPi=(\varPi_0,\varPi_1,\varPi_2)$ where $\varPi^i$ is the
posterior probability process of the non-zero drift being in the
spatial coordinate $i$ given $X$ for $i=0,1,2$. Due to $\sum_{i=0}^2
\varPi^i = 1$ clearly only two coordinates of $\varPi$ matter and
this is utilised by passing to the posterior probability ratio
process process $\varPhi=(\varPhi^1,\varPhi^2)$ defined by
$\varPhi^i = \varPi^i/\varPi^0$ for $i=1,2$. The processes $\varPi$
and $\varPhi$ stand in one-to-one correspondence and we study the
optimal stopping problem in terms of $\varPhi$ throughout. The
previous considerations take place under the probability measure
$\PP_{\!\pi} = \sum_{i=0}^2 \pi_i\;\! \PP_{\!i}$ where $\pi_i$ is
the prior probability of the non-zero drift being in the spatial
coordinate $i$ for $i=0,1,2$. In Section 4 we show that a measure
change from $\PP_{\!\pi}$ to $\PP_{\!0}$ simplifies the setting upon
verifying that the posterior probability ratio process $\varPhi^i$
coincides (up to the initial point) with the likelihood ratio
process $L^i$ of $\PP_{\!i}$ and $\PP_{\!0}$ given $X$ for $i=1,2$.
This provides an explicit link between the process $\varPhi$ and the
observed process $X$.

In Section 5 we show that the process $\varPhi$ solves a coupled
system of linear stochastic differential equations (of the geometric
Brownian motion type) driven by two independent Brownian motions.
This enables us to conclude that $\varPhi$ is a Markov/diffusion
process and derive a closed form expression for its infinitesimal
generator which is a second-order partial differential equation of
elliptic type. The optimal stopping problem for $\varPhi$ is Bolza
formulated and in Section 6 we disclose its Lagrange and Mayer
formulations (see \cite[Section 6]{PS} for the terminology). The
Lagrange formulation is expressed in terms of the local time of
$\varPhi$ on three straight lines which makes the optimal stopping
problem more intuitive.

The observed process $X$ is three-dimensional and in Section 7 we
consider the same optimal stopping problem when $X$ is
two-dimensional. In this case $\varPhi$ is a one-dimensional Markov/
diffusion process so that standard arguments enable us to solve the
optimal stopping problem in a closed form. The reduction of
dimension three to dimension two corresponds to either $\varPhi^1$
or $\varPhi^2$ becoming $0$ which is a natural boundary point for
both processes (cf.\ \cite{Fe}). The one-dimensional results of
Section 7 are used in Section 8 to derive existence of the optimal
stopping set and derive basic properties of the value function. We
show that the optimal stopping set consists of three convex sets
separated by the three straight lines that support the local time of
$\varPhi$ in the Lagrange formulation of the optimal stopping
problem. Using symmetry arguments combined with the one-dimensional
results of Section 7 we also derive the asymptotic behaviour of the
optimal stopping boundaries at zero and infinity.

In Section 9 we derive a directional smooth fit between the value
function and the loss function at the optimal stopping boundary. The
proof of the smooth fit makes use of the asymptotic behaviour of the
optimal stopping boundary at infinity to counter-balance the lack of
the global smoothness of the underlying loss function in the optimal
stopping problem. In Section 10 we show that the optimal stopping
boundaries are non-monotone in either direction of the state space
of $\varPhi$ and prove the existence of a `belly' which determines
their curvature/shape. These arguments rely on the general hint from
\cite[Remark 13]{Pe-3} on establishing the absence of jumps of the
optimal stopping boundaries and make use of Hopf's boundary point
lemma to derive a contradiction with the directional smooth fit.

In Section 11 we disclose the free-boundary problem which stands in
one-to-one correspondence with the optimal stopping problem and
establish the fact that the value function and the optimal stopping
boundaries solve the free-boundary problem uniquely. In Section 12
we show that the optimal stopping boundaries can be characterised as
the unique solution to a coupled system of nonlinear Fredholm
integral equations. These equations can be used to find the optimal
stopping boundaries numerically (using Picard iteration).

\section{Formulation of the problem}
%%%%%%%%%%%%%%%%%%%%%%%%%%%%%%%%%%%%

In this section we formulate the sequential testing problem under
consideration. The initial formulation of the problem will be
revaluated under a change of measure in the next section.

\vspace{6pt}

1.\ We consider a Bayesian formulation of the problem where it is
assumed that one observes a sample path of the three-dimensional
Brownian motion $X=(X^0,X^1,X^2)$, whose two coordinates $X^j$ and
$X^k$ are standard Brownian motions with zero drift, and the
remaining (uknown) coordinate $X^i$ is a standard Brownian motion
having a non-zero drift $\mu$ with a probability $\pi_i \in [0,1]$
for $i=0,1,2$ where $\pi_0 \p \pi_1 \p \pi_2 = 1$ and $i \ne j \ne
k$ belong to $\{0,1,2\}$. The problem is to detect which coordinate
is drifting as soon as possible and with minimal probabilities of
the wrong terminal decisions. This real-time detection
problem belongs to the class of sequential testing problems as
discussed in Section 1 above.

\vspace{6pt}

2.\ Standard arguments imply that the previous setting can be
realised on a probability space  $(\Omega,\cF,\PP_{\!\pi})$ with the
probability measure $\PP_{\!\pi}$ decomposed as follows
\begin{equation} \h{8pc} \label{3.1}
\PP_{\!\pi} = \pi_0\:\! \PP_{\!0} + \pi_1\:\! \PP_{\!1} + \pi_2\:\!
\PP_{\!2}
\end{equation}
for $\pi = (\pi_0,\pi_1,\pi_2) \in [0,1]^3$ satisfying $\pi_0 +
\pi_1 + \pi_2 = 1$ where $\PP_{\!i}$ is the probability measure
under which the observed process $X$ has the $i$\!-th coordinate
equal to a standard Brownian motion with drift $\mu$, and the
remaining two coordinates are standard Brownian motions with zero
drift for $i=0,1,2$, with the three coordinates being independent.
This can be formally achieved by introducing an unobservable random
variable $\theta$ taking values $0,1,2$ with probabilities
$\pi_0,\pi_1,\pi_2$ in $[0,1]$ satisfying $\pi_0 \p \pi_1 \p \pi_2 =
1$ and being independent from three standard Brownian motions
$B^0,B^1,B^2$ so that $X=(X^0,X^1,X^2)$ after starting at a point in
$\R^3$ solves the system of stochastic differential equations
\begin{equation} \h{8pc} \label{3.2}
dX_t^i = \mu\;\! I(\theta=i)\, dt + dB_t^i
\end{equation}
for $i=0,1,2$. Due to stationary and independent increments of
Brownian motion it is clear that the starting point of $X$ plays no
role in the sequel so we will leave it unspecified.

\vspace{6pt}

3.\ Being based upon the continued observation of $X$, the problem
is to test sequentially the hypotheses $H_0\, :\, \theta=0$, $H_1\,
:\, \theta=1$, $H_2\, :\, \theta=2$ with a minimal loss. For this,
we are given a sequential decision rule $(\tau,d_\tau)$, where
$\tau$ is a stopping time of $X$ (i.e.\ a stopping time with respect
to the natural filtration $\cF_t^X = \sigma(X_s\, \vert\, 0 \le s
\le t)$ of $X$ for $t \ge 0$), and $d_\tau$ is an
$\cF_\tau^X$\!-measurable random variable taking values in the set
$\{0,1,2\}$. After stopping the observation of $X$ at time $\tau$,
the terminal decision function $d_\tau$ takes value $i$ if and only
if the hypothesis $H_i$ is to be accepted for $i=0,1,2$. With a
constant $c>0$ given and fixed, the problem then becomes to compute
the risk function
\begin{equation} \h{0pc} \label{3.3}
V(\pi) = \inf_{(\tau,d_\tau)} \EE_\pi \Big[ \tau + c\;\! \Big(
I( \theta = 0 \:\! ,\:\! d_\tau \ne 0) \p I(\theta = 1 \:\! ,\:\!
d_\tau \ne 1) \p I(\theta = 2 \:\! ,\:\! d_\tau \ne 2) \Big) \Big]
\end{equation}
for $\pi = (\pi_0,\pi_1,\pi_2) \in [0,1]^3$ with $\pi_0 \p \pi_1 \p
\pi_2 = 1$ and find the optimal decision rule
$(\tau_*,d^*_{\tau_*})$ at which the infimum in \eqref{3.3} is
attained. Note that $\EE_\pi(\tau)$ in \eqref{3.3} is the expected
waiting time until the terminal decision is made, and
$\PP_{\!\pi}(\theta=i\:\! ,\:\! d_\tau \ne i)$ are probabilities of
the wrong terminal decisions for $i=0,1,2$. Clearly, each
probability $\PP_{\!\pi}(\theta=i\:\! ,\:\! d_\tau \ne i)$ could be
further decomposed into the sum of two probabilities
$\PP_{\!\pi}(\theta=i\:\! ,\:\! d_\tau = j)$ and
$\PP_{\!\pi}(\theta=i\:\! ,\:\! d_\tau = k)$ for $i=0,1,2$ and $i
\ne j \ne k$ in $\{0,1,2\}$, and each of the six resulting
probabilities could have a different constant/weight placed in front
of them, however, since the constrained problems are not considered
in the present paper as explained in the introduction, we only focus
on the canonical setting of a single constant/weight $c$ given in
\eqref{3.3} above.

\vspace{6pt}

4.\ To tackle the sequential testing problem \eqref{3.3} we consider
the \emph{posterior probability process} $\varPi = ((\varPi_t^0,
\varPi_t^1,\varPi_t^2))_{t \ge 0}$ of $H = (H_0,H_1,H_2)$ given $X$
that is defined by
\begin{equation} \h{9pc} \label{3.4}
\varPi_t^i = \PP_{\!\pi}(\theta = i\, \vert\, \cF_t^X)
\end{equation}
for $i=0,1,2$ and $t \ge 0$. Noting that for any decision rule
$(\tau,d_\tau)$ we have
\begin{equation} \h{1pc} \label{3.5}
\sum_{i=0}^2 \PP_{\!\pi}(\theta = i \:\! ,\:\! d_\tau \ne i) =
\sum_{i=0}^2 \EE_\pi \big[ \varPi_\tau^i\:\! I(d_\tau\! \ne\! i)
\big] = \sum_{i=0}^2 \EE_\pi \big[ (1 \m \varPi_\tau^i) \:\!
I(d_\tau\! =\! i) \big]
\end{equation}
where in the final equality we use that $\varPi_\tau^0 \p
\varPi_\tau^1 \p \varPi_\tau^2 = 1$, it follows that
\begin{align} \h{2pc} \label{3.6}
&\EE_\pi \Big[ \tau + c\;\! \Big( I( \theta = 0 \:\! ,\:\! d_\tau
\ne 0) \p I(\theta = 1 \:\! ,\:\! d_\tau \ne 1) \p I(\theta =
2 \:\! ,\:\! d_\tau \ne 2) \Big) \Big] \\[2pt] \notag &\h{2.6pc}
\ge \EE_\pi \Big[ \tau + c\;\! \Big( (1 \m \varPi_\tau^0) \wedge
(1 \m \varPi_\tau^1) \wedge (1 \m \varPi_\tau^2) \Big) \Big]
\end{align}
where equality is attained at the decision rule $(\tau,\tilde
d_\tau)$ with $\tilde d_\tau$ defined as follows
\begin{align} \h{5pc} \label{3.7}
\tilde d_\tau &= 0\;\;\; \text{if}\;\; (1 \m \varPi_\tau^0) \le
(1 \m \varPi_\tau^1) \wedge (1 \m \varPi_\tau^2) \\ \notag &=
1\;\;\; \text{if}\;\; (1 \m \varPi_\tau^1) \le (1 \m \varPi_
\tau^0) \wedge (1 \m \varPi_\tau^2) \\ \notag &= 2\;\;\; \text
{if}\;\; (1 \m \varPi_\tau^2) \le (1 \m \varPi_\tau^0) \wedge
(1 \m \varPi_\tau^1)\, .
\end{align}
This shows that the problem \eqref{3.3} is equivalent to the optimal
stopping problem
\begin{equation} \h{7pc} \label{3.8}
V(\pi) = \inf_\tau\;\! \EE_\pi \big[\;\! \tau + M(\varPi_\tau)
\;\! \big]
\end{equation}
where the infimum is taken over all stopping times $\tau$ of $X$, and
the function $M$ is given by
\begin{equation} \h{6pc} \label{3.9}
M(\pi) = c\;\! \big( (1 \m \pi_0)\! \wedge\! (1 \m \pi_1)\! \wedge\!
(1 \m \pi_2) \big)
\end{equation}
for $\pi = (\pi_0,\pi_1,\pi_2) \in [0,1]^3$ with $\pi_0 \p \pi_1 \p
\pi_2 = 1$. For this reason we focus on solving the optimal stopping
problem \eqref{3.8} in what follows.

\section{Measure change}
%%%%%%%%%%%%%%%%%%%%%%%%

In this section we show that changing the probability measure
$\PP_{\!\pi}$ for $\pi \in [0,1]^3$ with $\pi_0 +$ $ \pi_1 \p \pi_2
= 1$ to $\PP_{\!0}$ provides important simplifications of the
setting which make the subsequent analysis more transparent. The
change of measure argument is presented in Lemma 1 below. This is
then followed by a reformulation of the optimal stopping problem
\eqref{3.8} under the new probability measure $\PP_{\!0}$ in
Proposition 2 below.

\vspace{6pt}

1.\ To connect the process $\varPi$ in \eqref{3.8} to the observed
process $X$ we consider the \emph{likelihood ratio process} $L =
((L_t^1,L_t^2))_{t \ge 0}$ defined by
\begin{equation} \h{9pc} \label{4.1}
L_t^i = \frac{d \PP_{\!i,t}}{d \PP_{\!0,t}}
\end{equation}
where $\PP_{\!i,t}$ and $\PP_{\!0,t}$ denote the restrictions of
$\PP_{\!i}$ and $\PP_{\!0}$ to $\cF_t^X$ for $t \ge 0$ and $i=1,2$.
By the Girsanov theorem one finds that
\begin{equation} \h{9pc} \label{4.2}
L_t^i = e^{\;\!\mu\:\! (X_t^i - X_t^0)}
\end{equation}
for $t \ge 0$ and $i=1,2$. A direct calculation indicated below
shows that the \emph{posterior probability ratio process} $\varPhi =
((\varPhi_t^1,\varPhi_t^2))_{t \ge 0}$ defined by
\begin{equation} \h{9pc} \label{4.3}
\varPhi_t^i = \frac{\varPi_t^i}{\varPi_t^0}
\end{equation}
can be expressed in terms of $L$ (and hence $X$ as well) as follows
\begin{equation} \h{9pc} \label{4.4}
\varPhi_t^i = \varPhi_0^i\;\! L_t^i
\end{equation}
for $t \ge 0$ where $\varPhi_0^i = \pi_i/\pi_0$ for $i=1,2$.
Recalling that $\varPi_t^0 \p \varPi_t^1 \p \varPi_t^2 = 1$ and
formally setting $\varPhi_t^0 \equiv 1$ it is easily seen that
\eqref{4.3} is equivalent to
\begin{equation} \h{7.6pc} \label{4.5}
\varPi_t^i = \frac{\varPhi_t^i}{1 + \varPhi_t^1 + \varPhi_t^2}
\end{equation}
for $t \ge 0$ and $i=0,1,2$.

\vspace{6pt}

2.\ To derive \eqref{4.3}-\eqref{4.5} one may use a standard rule
for the Radon-Nikodym derivatives based on \eqref{3.1} that gives
\begin{align} \h{0pc} \label{4.6}
&\h{-0.3pc}\varPi_t^0 = \PP_{\!\pi}(\theta = 0\, \vert\, \cF_t^X) =
\sum_{i=0}^2 \pi_i\, \PP_{\!i}(\theta = 0\, \vert\, \cF_t^X)\, \frac
{d \PP_{\!i,t}}{d \PP_{\!\pi,t}} = \pi_0\;\! \frac{d \PP_{\!0,t}}{d
\PP_{\!\pi,t}} = \frac{1}{1 + \frac{\pi_1}{\pi_0} \frac{d \PP_{\!1,t}}
{d \PP_{\!0,t}} + \frac{\pi_2}{\pi_0} \frac{d \PP_{\!2,t}}{d \PP_{\!0,t}}}
\\ \label{4.7} &\h{-0.3pc}\varPi_t^1 = \PP_{\!\pi}(\theta = 1\, \vert\,
\cF_t^X) = \sum_{i=0}^2 \pi_i\, \PP_{\!i}(\theta = 1\, \vert\, \cF_t^X)
\, \frac{d \PP_{\!i,t}}{d \PP_{\!\pi,t}} = \pi_1\;\! \frac{d \PP_{\!1,t}}
{d \PP_{\!\pi,t}} = \frac{\frac{\pi_1}{\pi_0} \frac{d \PP_{\!1,t}}{d
\PP_{\!0,t}}}{1 + \frac{\pi_1}{\pi_0} \frac{d \PP_{\!1,t}}{d \PP_{\!0,t}}
+ \frac{\pi_2}{\pi_0} \frac{d \PP_{\!2,t}}{d \PP_{\!0,t}}} \\ \label{4.8}
&\h{-0.3pc}\varPi_t^2 = \PP_{\!\pi}(\theta = 2\, \vert\, \cF_t^X) =
\sum_{i=0}^2 \pi_i\, \PP_{\!i}(\theta = 2\, \vert\, \cF_t^X)\, \frac{d
\PP_{\!i,t}}{d \PP_{\!\pi,t}} = \pi_2\;\! \frac{d \PP_{\!2,t}}{d
\PP_{\!\pi,t}} = \frac{\frac{\pi_2}{\pi_0} \frac{d \PP_{\!2,t}}{d
\PP_{\!0,t}}}{1 + \frac{\pi_1}{\pi_0} \frac{d \PP_{\!1,t}}{d
\PP_{\!0,t}} + \frac{\pi_2}{\pi_0} \frac{d \PP_{\!2,t}}{d
\PP_{\!0,t}}}
\end{align}
where $\PP_{\!\pi,t}$ denotes the restriction of $\PP_{\!\pi}$ to
$\cF_t^X$ for $\pi = (\pi_0,\pi_1,\pi_2) \in [0,1]^3$ with $\pi_0 \p
\pi_1 \p$ $\pi_2 = 1$ and $t \ge 0$. It is then easily verified that
\eqref{4.6}-\eqref{4.8} imply \eqref{4.3}-\eqref{4.5} as claimed.

\vspace{6pt}

3.\ Previous arguments suggest that changing the probability measure
$\PP_{\!\pi}$ to $\PP_{\!0}$ appears to be of canonical interest in
the optimal stopping problem \eqref{3.8}. In the sequel we let
$\PP_{\!\pi,\tau}$ denote the restriction of $\PP_{\!\pi}$ to
$\cF_\tau^X$ where $\tau$ is a stopping time of $X$.

\vspace{12pt}

\textbf{Lemma 1.} \emph{The following identity holds
\begin{equation} \h{9pc} \label{4.9}
\frac{d \PP_{\!\pi,\tau}}{d \PP_{\!0,\tau}} = \frac{\pi_0}{
\varPi_\tau^0}
\end{equation}
for all stopping times $\tau$ of $X$ and all $\pi =
(\pi_0,\pi_1,\pi_2) \in [0,1]^3$ with $\pi_0 \p \pi_1 \p \pi_2 =
1$.}

\vspace{12pt}

\textbf{Proof.} Using the same arguments as in \eqref{4.6} above we
find that
\begin{equation} \h{2pc} \label{4.10}
\varPi_\tau^0 = \PP_{\!\pi}(\theta = 0\, \vert\, \cF_\tau^X) =
\sum_{i=0}^2 \pi_i\, \PP_{\!i}(\theta = 0\, \vert\, \cF_\tau^X)\,
\frac {d \PP_{\!i,\tau}}{d \PP_{\!\pi,\tau}} = \pi_0\;\! \frac{d
\PP_{\!0,\tau}}{d \PP_{\!\pi,\tau}}
\end{equation}
for any $\tau$ and $\pi$ as above. From \eqref{4.10} we see that
\eqref{4.9} holds and the proof is complete. \hfill $\square$

\vspace{6pt}

4.\ We now show that the optimal stopping problem \eqref{3.8} admits
a transparent reformulation under the probability measure
$\PP_{\!0}$ in terms of the process $\varPhi = (\varPhi^1,
\varPhi^2)$ defined in \eqref{4.3} above. Recall that $\varPhi^i$
starts at $\pi_i/\pi_0$ and this dependence on the initial point
will be indicated by a superscript $\pi_i/\pi_0$ to $\varPhi$
replacing its coordinate superscript $i$ for $i=1,2$ when needed.

\vspace{12pt}

\textbf{Proposition 2.} \emph{The value function $V$ from \eqref{3.8}
satisfies the identity
\begin{equation} \h{7pc} \label{4.11}
V(\pi) = \pi_0\;\! \hat V \Big( \frac{\pi_1}{\pi_0},\frac{\pi_2}
{\pi_0} \Big)
\end{equation}
where the value function $\hat V$ is given by
\begin{equation} \h{1pc} \label{4.12}
\hat V \Big( \frac{\pi_1}{\pi_0},\frac{\pi_2}{\pi_0} \Big) =
\inf_\tau\;\! \EE_0 \Big[ \int_0^\tau \big( 1 \p \varPhi_t^{
\pi_1/\pi_0}\! \p \varPhi_t^{\pi_2/\pi_0} \big)\, dt + \hat M
\big(\varPhi_\tau^{\pi_1/\pi_0}, \varPhi_\tau^{\pi_2/\pi_0}
\big)\:\! \Big]
\end{equation}
for $\pi = (\pi_0,\pi_1,\pi_2) \in [0,1]^3$ with $\pi_0 \p \pi_1 \p
\pi_2 = 1$ where
\begin{equation} \h{4pc} \label{4.13}
\hat M(\varphi_1,\varphi_2) = c\;\! \big( (\varphi_1 \p \varphi_2)
\wedge (1 \p \varphi_1) \wedge (1 \p \varphi_2) \big)
\end{equation}
for $(\varphi_1,\varphi_2) \in [0,\infty)^2$ and the infimum in
\eqref{4.12} is taken over all stopping times $\tau$ of $X$.}

\vspace{12pt}

\textbf{Proof.} For $\pi = (\pi_0,\pi_1,\pi_2) \in [0,1]^3$ with
$\pi_0 \p \pi_1 \p \pi_2 = 1$ given and fixed, it is enough to show
that the following identity holds
\begin{equation} \h{1pc} \label{4.14}
\EE_\pi \big[\;\! \tau + M(\varPi_\tau)\;\! \big] = \pi_0\,
\EE_0 \Big[ \int_0^\tau \big( 1 \p \varPhi_t^{ \pi_1/\pi_0}
\! \p \varPhi_t^{\pi_2/\pi_0} \big)\, dt + \hat M \big(\varPhi
_\tau^{\pi_1/\pi_0}, \varPhi_\tau^{\pi_2/\pi_0} \big)\:\! \Big]
\end{equation}
for all bounded stopping times $\tau$ of $X$. For this, suppose that
such a stopping time $\tau$ is given and fixed, and note by
\eqref{4.5}-\eqref{4.9} that
\begin{align} \h{1pc} \label{4.15}
\EE_\pi \big[\;\! \tau + M(\varPi_\tau)\;\! \big] &= \pi_0\, \EE_0
\Big[\,\frac{\tau}{\varPi_\tau^0} + \frac{M(\varPi_\tau)}{\Pi_\tau^0}
\, \Big] \\ \notag &= \pi_0\, \EE_0 \Big[ \tau \big( 1 \p \varPhi_
\tau^{ \pi_1/\pi_0}\! \p \varPhi_\tau^{\pi_2/\pi_0} \big) + \hat M
\big(\varPhi _\tau^{\pi_1/\pi_0}, \varPhi_\tau^{\pi_2/\pi_0} \big)
\:\! \Big]\, .
\end{align}
Setting $M_t = 1 \p \varPhi_t^{ \pi_1/\pi_0}\! \p
\varPhi_t^{\pi_2/\pi_0}$ for $t \ge 0$ we see by \eqref{4.2} and
\eqref{4.4} that $M = (M_t)_{t \ge 0}$ is a continuous martingale
under $\PP_{\!0}$ so that integration by parts gives
\begin{equation} \h{7pc} \label{4.16}
\textstyle t M_t = \int_0^t M_s\, ds + \int_0^t s \, dM_s
\end{equation}
where the final term defines a continuous martingale under
$\PP_{\!0}$ for $t \ge 0$. By the optional sampling theorem we
therefore get
\begin{equation} \h{7pc} \label{4.17}
\textstyle \EE_0 \big( \tau M_\tau \big) = \EE_0 \big(\int_0^\tau
M_t\, dt\:\! \big)\, .
\end{equation}
Inserting this back into \eqref{4.15} we obtain \eqref{4.14} as
claimed and the proof is complete. \hfill $\square$

\vspace{12pt}

5.\ It is clear from \eqref{4.2} and \eqref{4.4} that $\varPhi =
(\varPhi^1,\varPhi^2)$ is a strong Markov/diffusion process. We will
formally verify this fact in the next section by deriving a coupled
system of stochastic differential equations (driven by two
independent Brownian motions) that $\varPhi$ solves. Denoting the
probability law of $\varPhi^\varphi = (\varPhi^{\varphi_1},
\varPhi^{\varphi_2})$ under $\PP_{\!0}$ by $\PP_{\!\varphi}^0 =
\PP_{\!\varphi_1,\varphi_2}^0$ (where we move $0$ from the subscript
to a superscript for notational reasons) we see that the optimal
stopping problem \eqref{4.12} can be rewritten as follows
\begin{equation} \h{1pc} \label{4.18}
\hat V (\varphi_1,\varphi_2) = \inf_\tau\;\! \EE_{\varphi_1,
\varphi_2}^0 \Big[ \int_0^\tau\! \big( 1 \p \varPhi_t^1 \p
\varPhi_t^2 \big)\, dt + c\:\! \Big( (\varPhi_\tau^1 \p
\varPhi_\tau^2) \wedge (1 \p \varPhi_\tau^1) \wedge (1 \p
\varPhi_\tau^2) \Big)\:\! \Big]
\end{equation}
for $(\varphi_1,\varphi_2) \in [0,\infty)^2$ with
$\PP_{\!\varphi_1,\varphi_2} \big( (\varPhi_0^1, \varPhi_0^2)\! =\!
(\varphi_1,\varphi_2) \big) = 1$ where the infimum in \eqref{4.18}
is taken over all stopping times $\tau$ of $\varPhi$. In this way we
have reduced the initial sequential testing problem \eqref{3.3} to
the optimal stopping problem \eqref{4.18} for the strong
Markov/diffusion process $\varPhi$. We will see in the next section
that this optimal stopping problem is inherently/fully
two-dimensional with the infinitesimal generator of $\varPhi$ being
of elliptic type.

\section{Elliptic PDE}
%%%%%%%%%%%%%%%%%%%%%%

In this section we derive a coupled system of stochastic
differential equations (driven by two independent Brownian motions)
that $\varPhi = (\varPhi^1,\varPhi^2)$ solves. From this system we
derive a closed-form expression for the infinitesimal generator of
$\varPhi$ that can be recognised as a partial differential equations
of elliptic type. We also show that a diffeomorphic transformation
of logarithmic type maps the process $\varPhi$ (\:\!and its state
space $(0,\infty)^2$\!) to a process $Z$ (\:\!and its state space
$\R^2$\!) whose coordinate processes $Z^1$ and $Z^2$ are independent
Brownian motions with a non-zero and zero drift respectively.

\vspace{6pt}

1.\ From \eqref{4.2} and \eqref{4.4} we see that
\begin{equation} \h{3pc} \label{5.1}
\varPhi_t^1 = \varphi_1\:\! e^{\;\! \mu\:\! (B_t^1 - B_t^0) - \mu^2
t} \;\;\; \&\;\;\; \varPhi_t^2 = \varphi_2\:\! e^{\;\! \mu\:\! (B_t^2
- B_t^0) - \mu^2 t}
\end{equation}
under $\PP_{\!0}$ for $t \ge 0$ where $\varphi_1$ and $\varphi_2$
belong to $[0,\infty)$. Hence by It\^o's formula we find that
\begin{align} \h{7pc} \label{5.2}
&d \varPhi_t^1 = \mu\:\! \varPhi_t^1 (dB_t^1 \m dB_t^0) \\[4pt]
\label{5.3} &d \varPhi_t^2 = \mu\:\! \varPhi_t^2 (dB_t^2 \m dB_t^0)
\end{align}
under $\PP_{\!0}$ with $\varPhi_0^1 = \varphi_1$ and $\varPhi_0^2 =
\varphi_2$ in $[0,\infty)$. This shows that  $\varPhi^1$ and
$\varPhi^2$ are two correlated geometric Brownian motions.

\vspace{6pt}

2.\ A well-known (and easily verifiable) fact states that if $\tilde
B^1$ and $\tilde B^2$ are two correlated standard Brownian motions
satisfying $\EE(\tilde B_t^1 \tilde B_t^2) = \rho\:\! t$ for $t \ge
0$ with $\rho \in (-1,1)$, then $(\tilde B^1 +$ $\tilde B^2)/\sqrt{2
(1 \p \rho)}$ and $(\tilde B^1 \m \tilde B^2)/\sqrt{2 (1 \m \rho)}$
are two independent standard Brownian motions. Applying this
implication to $\tilde B^1 := (B^1 \m B^0)/\sqrt{2}$ and $\tilde B^2
:= (B^2 \m B^0)/\sqrt{2}$ with $\rho = 1/2$ it follows that
\begin{equation} \h{1pc} \label{5.4}
W^1 := \frac{\tilde B^1 \p \tilde B^2}{\sqrt{3}} = \frac{B^1 \p B^2
\m 2 B^0}{\sqrt{6}} \quad \& \quad W^2 := \frac{\tilde B^1 \m \tilde
B^2}{1} = \frac{B^1 \m B^2}{\sqrt{2}}
\end{equation}
are two independent standard Brownian motions. From \eqref{5.4} we
see that
\begin{equation} \h{1pc} \label{5.5}
\tilde B^1 := \frac{B^1 \m B^0}{\sqrt{2}} = \frac{\sqrt{3}\;\! W^1 \p
W^2}{2} \quad \& \quad \tilde B^2 := \frac{B^2 \m B^0}{\sqrt{2}}
= \frac{\sqrt{3}\;\! W^1 \m W^2}{2}\, .
\end{equation}

\vspace{6pt}

3.\ Making use of \eqref{5.5} in \eqref{5.2}+\eqref{5.3} we obtain
\begin{align} \h{5.5pc} \label{5.6}
&d \varPhi_t^1 = \frac{\mu}{\sqrt{2}}\, \varPhi_t^1\:\! \big(\sqrt{3}
\;\! dW_t^1 \p dW_t^2 \big) \\[1pt] \label{5.7} &d \varPhi_t^2 =
\frac{\mu}{\sqrt{2}}\, \varPhi_t^2\:\! \big(\sqrt{3}\;\! dW_t^1
\m dW_t^2 \big)
\end{align}
with $\varPhi_0^1 = \varphi_1$ and $\varPhi_0^2 = \varphi_2$ in
$[0,\infty)$. This is a coupled system of stochastic differential
equations (\:\!driven by two independent standard Brownian motions
$W^1$ and $W^2$\!) that $\varPhi^1$ and $\varPhi^2$ solve (strongly)
and this solution is pathwise unique (see e.g.\ \cite[pp
128-131]{RW}). Moreover, the solution $\varPhi =
(\varPhi^1,\varPhi^2)$ is both a strong Markov process (see e.g.\
\cite[pp 158-163]{RW}) and a strong Feller process (see e.g.\
\cite[pp 170-173]{RW}). Making use of \eqref{5.5} in \eqref{5.1} we
see that
\begin{equation} \h{3pc} \label{5.8}
\varPhi_t^1 = \varphi_1\:\! e^{\;\! \frac{\mu}{\sqrt{2}}\:\! (\sqrt{3}
\:\! W_t^1 + W_t^2) - \mu^2 t} \;\;\; \&\;\;\; \varPhi_t^2 = \varphi_2
\:\! e^{\;\! \frac{\mu}{\sqrt{2}}\:\! (\sqrt{3} \:\! W_t^1 - W_t^2)
- \mu^2 t}
\end{equation}
under $\PP_{\!0}$ for $t \ge 0$ where $\varphi_1$ and $\varphi_2$
belong to $[0,\infty)$. Often we will write $\varPhi_t^{\varphi_1}$
and $\varPhi_t^{\varphi_2}$ for $t \ge 0$ to indicate dependence of
$\varPhi^1$ and $\varPhi^2$ on the initial points $\varphi_1$ and
$\varphi_2$ in $[0,\infty)$.

\vspace{6pt}

4.\ Knowing that $\varPhi = (\varPhi^1,\varPhi^2)$ solves the system
\eqref{5.6}+\eqref{5.7} and making use of It\^o's calculus we find
that the infinitesimal generator of $\varPhi$ is given by
\begin{equation} \h{4.5pc} \label{5.9}
\LL_\varPhi = \mu^2 \big(\;\! \varphi_1^2\;\! \partial_{\varphi_1
\varphi_1}^2\! + \varphi_1 \varphi_2\;\! \partial_{\varphi_1 \varphi_2}
^2\! + \varphi_2^2\;\! \partial_{\varphi_2 \varphi_2}^2\;\! \big)
\end{equation}
for $\varphi_1$ and $\varphi_2$ in $(0,\infty)$ (see e.g.\ (2.7) in
\cite{Pe-3}). A standard classification of partial differential
equations shows that $\LL_\varPhi$ is of elliptic type (see e.g.\
(2.12) in \cite{Pe-3}).

\vspace{6pt}

5.\ Defining a diffeomorphic transformation of $(0,\infty)^2$ to
$\R^2$ by
\begin{equation} \h{4.5pc} \label{5.10}
D(\varphi_1,\varphi_2) = \big( \log(\varphi_1 \varphi_2), \log(
\varphi_1 / \varphi_2) \big)
\end{equation}
for $(\varphi_1,\varphi_2) \in (0,\infty)^2$, and setting
\begin{equation} \h{3pc} \label{5.11}
Z = (Z^1,Z^2) = D(\varPhi_1,\varPhi_2) = \big( \log(\varPhi_1
\varPhi_2), \log( \varPhi_1 / \varPhi_2) \big)
\end{equation}
we see from \eqref{5.8} that
\begin{equation} \h{3pc} \label{5.12}
Z_t^1 = Z_0^1 - 2 \mu^2 t + \sqrt{6}\:\! \mu\:\! W_t^1 \quad \& \quad
Z_t^2 = Z_0^2 + \sqrt{2}\:\! \mu\:\! W_t^2
\end{equation}
under $\PP_{\!0}$ for $t \ge 0$ with $Z_0^1 = \log(\varphi_1
\varphi_2)$ and $Z_0^2 = \log(\varphi_1 / \varphi_2)$. This
establishes a one-to-one correspondence between the process
$\varPhi$ in $(0,\infty)^2$ and the process $Z$ in $\R^2$. Although
the latter process $Z$ may be viewed as a canonical building block
which further clarifies the underlying setting, we will mainly study
the optimal stopping problem \eqref{4.18} by means of the former
process $\varPhi$ in the sequel.

\section{Lagrange and Mayer formulations}
%%%%%%%%%%%%%%%%%%%%%%%%%%%%%%%%%%%%%%%%%

The optimal stopping problem \eqref{4.18} is Bolza formulated. In
this section we derive its Lagrange and Mayer reformulations which
are helpful in the subsequent analysis of the problem.

\vspace{6pt}

1.\ We first consider the Lagrange reformulation of the optimal
stopping problem \eqref{4.18}. For this, note that the loss function
$\hat M$ from \eqref{4.13} that appears on the right-hand side of
\eqref{4.18} is not smooth at the three straight lines
\begin{align} \h{3pc} \label{6.1}
&c_0 = \{\, (\varphi_1,\varphi_2) \in [0,\infty)^2\; \vert\; \varphi_1
=1\;\; \&\;\; \varphi_2 \in [0,1]\,\} \\[3pt] \label{6.2} &c_1 = \{\,
(\varphi_1,\varphi_2) \in [0,\infty)^2\; \vert\; \varphi_1 \in [0,1]
\;\; \&\;\; \varphi_2 = 1\,\} \\[3pt] \label{6.3} &c_2 = \{\, (\varphi_1,
\varphi_2) \in [0,\infty)^2\; \vert\; \varphi_1 = \varphi_2 \in
[1,\infty)\,\}
\end{align}
ordered clockwise (see Figure 1 below). Note moreover that $\hat M$
is linear off the three straight lines and given by
\begin{align} \h{4pc} \label{6.4}
\hat M(\varphi_1,\varphi_2) &= c (\varphi_1 \p \varphi_2)\;\; \text{for}
\;\; (\varphi_1,\varphi_2) \in \Delta_0 \\[3pt] \notag &= c (1 \p \varphi_1)
\;\; \text{for} \;\; (\varphi_1,\varphi_2) \in \Delta_1 \\[2pt]\notag &= c
(1 \p \varphi_2)\;\; \text{for} \;\; (\varphi_1,\varphi_2) \in \Delta_2
\end{align}
where $\Delta_0 := [0,1]^2$ is a subset of the state space
surrounded by $c_0$ and $c_1$ (from the right and above), $\Delta_1
:= \{\, (\varphi_1,\varphi_2) \in [0,\infty)^2\; \vert\; \varphi_2
\ge \varphi_1 \ge 1\; \}$ is a subset of the state space surrounded
by $c_1$ and $c_2$ (from below and the right), and $\Delta_2 := \{\,
(\varphi_1,\varphi_2) \in [0,\infty)^2\; \vert\; \varphi_1 \ge
\varphi_2 \ge 1\; \}$ is a subset of the state space surrounded by
$c_0$ and $c_2$ (from the left and above).

\vspace{12pt}

\textbf{Proposition 3.} \emph{The value function $\hat V$ from
\eqref{4.18} can be expressed as
\begin{equation} \h{-0.5pc} \label{6.5}
\hat V(\varphi_1,\varphi_2) = \inf_\tau\;\! \EE_{\varphi_1,
\varphi_2}^0 \Big[ \int_0^\tau\! \big( 1 \p \varPhi_t^1 \p
\varPhi_t^2 \big)\, dt - \frac{c}{2}\:\! \Big( \ell_\tau^
{c_0}(\varPhi) \p \ell_\tau^{c_1}(\varPhi) \p \ell_\tau^
{c_2}(\varPhi)\Big)\:\! \Big] + \hat M(\varphi_1,\varphi_2)
\end{equation}
for $(\varphi_1,\varphi_2) \in [0,\infty)^2$ where
$\ell^{c_i}(\varPhi)$ is the local time of $\varPhi$ at $c_i$ for
$i=0,1,2$ given by
\begin{align} \h{0.5pc} \label{6.6}
&\ell_\tau^{c_0}(\varPhi) = \PP\:\! \text{-} \lim_{\eps \downarrow
0}\, \frac{1}{2 \eps} \int_0^\tau\! I(1 \m \eps < \varPhi_t^1 < 1
\p \eps)\, I( 0 \le \varPhi_t^2 \le 1)\: d \langle \varPhi^1,
\varPhi^1 \rangle_t \\[2pt] \label{6.7} &\ell_\tau^{c_1}(\varPhi)
= \PP\:\! \text{-} \lim_{\eps \downarrow 0}\, \frac{1}{2 \eps}
\int_0^\tau\! I(1 \m \eps < \varPhi_t^2 < 1 \p \eps)\, I( 0 \le
\varPhi_t^1 \le 1)\: d \langle \varPhi^2,\varPhi^2 \rangle_t
\\[2pt] \label{6.8} &\ell_\tau^{c_2}(\varPhi) = \PP\:\! \text{-}
\lim_{\eps \downarrow 0}\, \frac{1}{2 \eps} \int_0^\tau\! I(-\eps
< \varPhi_t^2 \m \varPhi_t^1 < \eps)\, I( \varPhi_t^1 \ge 1,
\varPhi_t^2 \ge 1)\: d \langle \varPhi^2 \m \varPhi^1,\varPhi^2
\m \varPhi^1 \rangle_t
\end{align}
and the infimum in \eqref{6.5} is taken over all stopping times
$\tau$ of $\varPhi$.}

\vspace{12pt}

\textbf{Proof.} It is evident from \eqref{6.4} that $\hat M$
restricted to $\Delta_0 \cup \Delta_2$ can be extended to a twice
continuously differentiable function $\hat F$ on $[0,\infty)^2\!
\setminus\! c_0$. Then $\hat M = \hat F + (\hat M \m \hat F)$ and
$\hat M^1 := \hat F$ is not smooth at $c_0$ while $\hat M^2 := \hat
M \m \hat F$ is not smooth at $c_1$ and $c_2$. Since $c_0$ is the
graph of a (linear) function of $\varphi_2$, and $c_1$ and $c_2$ are
the graphs of (linear) functions of $\varphi_1$, we see that the
change-of-variable formula with local time on surfaces \cite[Theorem
2.1]{Pe-2} is applicable to $\hat M^1$ and $\hat M^2$ composed with
$\varPhi$, \vspace{-1pt} where we note that $\hat
M_{\varphi_1}^1(\varphi_1+,\varphi_2) \m \hat
M_{\varphi_1}^1(\varphi_1-,\varphi_2) = -c$ for
$(\varphi_1,\varphi_2) \in c_0$ and $\hat
M_{\varphi_2}^2(\varphi_1,\varphi_2+) \m \hat
M_{\varphi_2}^2(\varphi_1,\varphi_2-) = -c$ for
$(\varphi_1,\varphi_2) \in c_1 \cup c_2$. Hence the formula is also
applicable to $\hat M$ composed with $\varPhi$ and this gives
\begin{align} \h{0.5pc} \label{6.9}
\hat M(\varPhi_t^1,\varPhi_t^2) &= \hat M(\varphi_1,\varphi_2) +
c \int_0^t\! I(\varPhi_s\! \in\! \Delta_0 \cup \Delta_1)\, d \varPhi_s^1 +
c \int_0^t\! I(\varPhi_s\! \in\! \Delta_0 \cup \Delta_2)\, d \varPhi_s^2 \\
\notag &\h{13pt}- \frac{c}{2}\:\! \big(\:\! \ell_t^{c_0}(\varPhi)
\p \ell_t^{c_1} (\varPhi) \p \ell_t^{c_2}(\varPhi)\:\! \big)
\end{align}
for $(\varphi_1,\varphi_2) \in [0,\infty)^2$ and $t \ge 0$ where the
local times are defined in \eqref{6.6}-\eqref{6.8} above. Since
$\varPhi^1$ and $\varPhi^2$ are continuous martingales under
$\PP_{\!0}$ we see that the two integrals on the right-hand side of
\eqref{6.9} are continuous martingales under $\PP_{\!0}$ as well. By
the optional sampling theorem we therefore find from \eqref{6.9}
that
\begin{equation} \h{2pc} \label{6.10}
\EE_{\varphi_1,\varphi_2}^0 \big[ \hat M(\varPhi_\tau^1,\varPhi_\tau^2)
\big] = \hat M(\varphi_1,\varphi_2) - \frac{c}{2}\, \EE_{\varphi_1,
\varphi_2}^0 \big[\:\! \ell_\tau^{c_0}(\varPhi) \p \ell_\tau^{c_1}
(\varPhi) \p \ell_\tau^{c_2}(\varPhi)\:\! \big]
\end{equation}
for all $(\varphi_1,\varphi_2) \in [0,\infty)^2$ and all stopping
times $\tau$ of $\varPhi$. Inserting \eqref{6.10} into \eqref{4.18}
we obtain \eqref{6.5} as claimed and the proof is complete. \hfill
$\square$

\vspace{12pt}

The Lagrange reformulation \eqref{6.5} of the optimal stopping
problem \eqref{4.18} reveals the underlying rationale for continuing
vs stopping in a clearer manner. Indeed, recalling that the local
time process $t \mapsto \ell_t^{c_i}(\varPhi)$ strictly increases
only when $\varPhi$ is at $c_i$, and that $\ell_t^{c_i}(\varPhi)
\sim \sqrt{t}$ is strictly larger than $\int_0^t (1 \p \varPhi_s^1
\p \varPhi_s^2)\, ds \sim t$ for small $t$, we see from \eqref{6.5}
that it should never be optimal to stop at $c_i$ and the incentive
for stopping should increase the further away $\varPhi$ gets from
$c_i$ for $=0,1,2$. We will see in Section 8 below that these
informal conjectures can be formalised and this will give a proof of
the fact that the three straight lines $c_0, c_1, c_2$ are contained
in the continuation set of the optimal stopping problem
\eqref{4.18}.

\vspace{12pt}

2.\ We next consider the Mayer reformulation of the optimal stopping
problem \eqref{4.18}. For this, in addition to $\hat M$ in
\eqref{4.13} above, define
\begin{equation} \h{2pc} \label{6.11}
\check M(\varphi_1,\varphi_2) = \frac{1}{\mu^2}\:\! \Big( \varphi_1
\big( \log \varphi_1 \m 1 \big) + \varphi_2 \big( \log \varphi_2 \m 1
\big) - \frac{1}{2} \log(\varphi_1 \varphi_2) \Big)
\end{equation}
and set $M(\varphi_1,\varphi_2) = \check M(\varphi_1,\varphi_2) +
\hat M(\varphi_1,\varphi_2)$ for $(\varphi_1,\varphi_2) \in
(0,\infty)^2$.

\vspace{12pt}

\textbf{Proposition 4.} \emph{The value function $\hat V$ from
\eqref{4.18} can be expressed as
\begin{equation} \h{4pc} \label{6.12}
\hat V(\varphi_1,\varphi_2) = \inf_\tau\;\! \EE_{\varphi_1,
\varphi_2}^0 \big[ M(\varPhi_\tau^1,\varPhi_\tau^2)\:\! \big]
- \check M(\varphi_1,\varphi_2)
\end{equation}
for $(\varphi_1,\varphi_2) \in (0,\infty)^2$ where the infimum is
taken over all stopping times $\tau$ of $\varPhi$.}

\vspace{12pt}

\textbf{Proof.} Recalling the closed-form expression for
$\LL_\varPhi$ in \eqref{5.9} it is easily verified that
\begin{equation} \h{6pc} \label{6.13}
\LL_\varPhi\:\! \check M(\varphi_1,\varphi_2) = 1 \p \varphi_1
\p \varphi_2
\end{equation}
for $(\varphi_1,\varphi_2) \in (0,\infty)^2$. By It\^o's formula we
thus find using \eqref{5.6}+\eqref{5.7} above that
\begin{align} \h{1pc} \label{6.14}
&\check M(\varPhi_t^1,\varPhi_t^2) = \check M(\varphi_1,\varphi_2) +
\int_0^t \check M_{\varphi_1}(\varPhi_s^1,\varPhi_s^2)\, d \varPhi_s^1
+ \int_0^t \check M_{\varphi_2}(\varPhi_s^1,\varPhi_s^2)\, d \varPhi
_s^2 \\ \notag &+ \int_0^t \LL_\varPhi\:\! \check M(\varPhi_s^1,
\varPhi_s^2) \, ds = \check M(\varphi_1,\varphi_2) + \int_0^t
\frac{\mu}{\sqrt{2}} \big( \varPhi_s^1 \log \varPhi_s^1 \m
\tfrac{1}{2} \big) \big( \sqrt{3}\, d W_s^1 \p d W_s^2 \big) \\
\notag &+ \int_0^t \frac{\mu}{\sqrt{2}} \big( \varPhi_s^2 \log
\varPhi_s^2 \m \tfrac{1}{2} \big) \big( \sqrt{3}\, d W_s^1 \m d
W_s^2 \big) + \int_0^t \big(1 \p \varPhi_s^1 \p \varPhi_s^2 \big)
\, ds
\end{align}
for $(\varphi_1,\varphi_2) \in (0,\infty)^2$ and $t \ge 0$ where the
two integrals on the right-hand side define continuous local
martingales under $\PP_{\!0}$. Making use of a localisation sequence
of stopping times for these two local martingales if needed, and
applying the optional sampling theorem, we find from \eqref{6.14}
that
\begin{equation} \h{2pc} \label{6.15}
\EE_{\varphi_1,\varphi_2}^0 \big[ \check M(\varPhi_\tau^1,\varPhi_
\tau^2) \big] = \check M(\varphi_1,\varphi_2) + \EE_{\varphi_1,
\varphi_2}^0 \Big[ \int_0^\tau\! \big( 1 \p \varPhi_t^1 \p
\varPhi_t^2 \big)\, dt\:\! \Big]
\end{equation}
for all $(\varphi_1,\varphi_2) \in (0,\infty)^2$ and all (bounded)
stopping times $\tau$ of $\varPhi$. Inserting \eqref{6.15} into
\eqref{4.18} we obtain \eqref{6.12} as claimed and the proof is
complete. \hfill $\square$

\section{Two dimensions}
%%%%%%%%%%%%%%%%%%%%%%%%

The observed process $X$ in the initial sequential testing problem
\eqref{3.3} is three-dimensional. In this section we consider the
analogue of \eqref{3.3} and the resulting optimal stopping problem
\eqref{4.18} when $X$ is two-dimensional. The reduction of dimension
three to dimension two corresponds to either $\varPhi^1$ or
$\varPhi^2$ becoming $0$ which is a natural boundary point for both
processes (cf.\ \cite{Fe}). This shows that $\varPhi$ is a
one-dimensional Markov/diffusion process when $X$ is two-dimensional
so that standard arguments enable us to solve the problem
\eqref{4.18} in a closed form. The derived results for the
one-dimensional optimal stopping problem \eqref{4.18} when $X$ is
two-dimensional will be used in the subsequent analysis of the
two-dimensional optimal stopping problem \eqref{4.18} when $X$ is
three-dimensional.

\vspace{6pt}

1.\ Using the same arguments as above, it is easily seen that the
sequential testing problem \eqref{3.3} when $X$ is two-dimensional
reduces to the optimal stopping problem \eqref{4.18} with
$\varPhi^2$ being formally equal to zero. Omitting the subscript $1$
from $\varphi_1$ for simplicity, we thus see that the optimal
stopping problem \eqref{4.18} reads
\begin{equation} \h{5pc} \label{7.1}
\hat V (\varphi) = \inf_\tau\;\! \EE_{\varphi}^0 \Big[ \int_0^\tau\!
\big( 1 \p \varPhi_t \big)\, dt + c\:\! \big( 1\! \wedge\! \varPhi_\tau
\big) \:\! \Big]
\end{equation}
for $\varphi \in [0, \infty)$ with $\PP_{\!\varphi}^0(\varPhi_0\!
=\! \varphi) = 1$ where the infimum in \eqref{7.1} is taken over all
stopping times $\tau$ of $\varPhi$. From \eqref{5.1} and \eqref{5.2}
we see that
\begin{align} \h{7pc} \label{7.2}
&\varPhi_t = \varphi \, e^{\sqrt{2}\:\! \mu\:\! W_t - \mu^2 t}
\\[2pt] \label{7.3} &d \varPhi_t = \sqrt{2}\;\! \mu\;\! \varPhi_t
\, dW_t
\end{align}
under $\PP_{\!0}$ for $t \ge 0$ with $\varPhi_0 = \varphi$ in
$[0,\infty)$ where $W := (B^1 \m B^0)/\sqrt{2}$ is a standard Brown-
ian motion. From \eqref{7.3} we see that the infinitesimal generator
of $\varPhi$ is given by
\begin{equation} \h{7pc} \label{7.4}
\LL_\varPhi = \mu^2\:\! \varphi^2 \frac{d^2}{d \varphi^2}
\end{equation}
which also follows formally by setting $\varphi_2 = 0$ in
\eqref{5.9} above.

Recognising the loss function in \eqref{7.1} as $\hat M(\varphi) =
c\:\! (1\! \wedge\! \varphi)$ for $\varphi \in [0,\infty)$, standard
arguments imply (see e.g.\ \cite{PS}) that $\hat V$ should solve the
free-boundary problem
\begin{align} \h{5pc} \label{7.5}
&\LL_\varPhi \hat V(\varphi) = -(1 \p \varphi)\;\; \text{for}\;\;
\varphi \in (\varphi_0^*,\varphi_1^*) \\[1pt] \label{7.6} &\hat
V(\varphi_i^*) = \hat M(\varphi_i^*)\;\; \text{for}\;\; i=0,1
\;\; \text{(instantaneous stopping)} \\[1pt] \label{7.7} &\hat
V'(\varphi_i^*) = \hat M'(\varphi_i^*)\;\; \text{for}\;\; i=0,
1\;\; \text{(smooth fit)}
\end{align}
where $0 < \varphi_0^* < 1 < \varphi_1^* < \infty$ are the optimal
stopping/boundary points to be found and we have $\hat V(\varphi) =
\hat M(\varphi)$ for $\varphi \in [0,\varphi_0^*) \cup
(\varphi_1^*,\infty)$ as well (in addition to \eqref{7.6} above).

The general solution to the ordinary differential equation
\eqref{7.5} is given by
\begin{equation} \h{6pc} \label{7.8}
\hat V(\varphi) = A\:\! \varphi + B + \frac{1}{\mu^2}\:\! (1 \m
\varphi) \log \varphi
\end{equation}
for $\varphi > 0$ where $A$ and $B$ are two undetermined real
constants. Boundary conditions \eqref{7.6} and \eqref{7.7} then read
as follows
\begin{align} \h{6pc} \label{7.9}
&A\:\! \varphi_0^* + B + \frac{1}{\mu^2}\:\! (1 \m \varphi_0^*) \log
\varphi_0^* = c\:\! \varphi_0^* \\ \label{7.10}
&A\:\! \varphi_1^* + B + \frac{1}{\mu^2}\:\! (1 \m \varphi_1^*) \log
\varphi_1^* = c \\ \label{7.11}
&A + \frac{1}{\mu^2}\:\! \Big( \frac{1}{\varphi_0^*} - \log \varphi_0^*
- 1 \big) = c \\ \label{7.12}
&A + \frac{1}{\mu^2}\:\! \Big( \frac{1}{\varphi_1^*} - \log \varphi_1^*
- 1 \big) = 0\, .
\end{align}
It is a matter of routine to verify that the system
\eqref{7.9}-\eqref{7.12} has a unique solution given by
\begin{equation} \h{2pc} \label{7.13}
A^* = c - \frac{1}{\mu^2}\:\! \Big( \frac{1}{\varphi_0^*} - \log
\varphi_0^* - 1 \Big)\;\;\; \&\;\;\; B^* = c - \frac{1}{\mu^2}\:\!
\Big( \varphi_1^* + \log \varphi_1^* - 1 \Big)\, .
\end{equation}
where $\varphi_0^*$ and $\varphi_1^*$ are the unique solution to
\begin{equation} \h{2pc} \label{7.14}
\frac{1}{\mu^2}\:\! \Big( \frac{1}{\varphi_0^*} - \frac{1}{\varphi_1^*}
+ \log \Big( \frac{\varphi_1^*}{\varphi_0^*} \Big) \Big) = c\;\;\;
\& \;\;\; \frac{1}{\mu^2}\:\! \Big( \varphi_1^* - \varphi_0^* +
\log \Big( \frac{\varphi_1^*}{\varphi_0^*} \Big) \Big) = c
\end{equation}
satisfying $0 < \varphi_0^* < 1 < \varphi_1^* < \infty$.

By symmetry we may conclude that $\varphi_0^* = 1/\varphi_1^*$ so
that \eqref{7.13} and \eqref{7.14} reduce to
\begin{align} \h{6pc} \label{7.15}
&A^* = B^* = c - \frac{1}{\mu^2}\:\! \Big( \varphi_1^* + \log \varphi_1^*
- 1 \Big) \\ \label{7.16} &\varphi_1^* - \frac{1}{\varphi_1^*} + 2\:\! \log
\varphi_1^* = c\:\! \mu^2
\end{align}
respectively. It follows from \eqref{7.8} and \eqref{7.15} that
\begin{align} \h{0pc} \label{7.17}
\hat V^*(\varphi) &= \Big( c - \frac{1}{\mu^2}\:\! \Big( \varphi_1^*
+ \log \varphi_1^* - 1 \Big) \Big) (1 \p \varphi) + \frac{1}{\mu^2}
\:\! (1 \m \varphi) \log \varphi\;\; \text{for}\;\; \varphi \in (1/
\varphi_1^*,\varphi_1^*) \\ \notag &= \hat M(\varphi)\;\; \text{for}
\;\; \varphi \in [0,1/\varphi_1^*] \cup [\varphi_1^*,\infty)
\end{align}
defines a candidate value function for the optimal stopping problem
\eqref{7.1}.

Applying the It\^o-Tanaka formula (cf.\ \cite[p.\ 223]{RY}) to $\hat
V^*$ composed with $\varPhi$, which reduces to It\^o's formula due
to smooth fit \eqref{7.7}, and making use of the optional sampling
theorem, it is easily verified that $\hat V^*$ from \eqref{7.17}
coincides with the value function $\hat V$ from \eqref{7.1} and the
optimal stopping time (at which the infimum in \eqref{7.1} is
attained) is given by
\begin{equation} \h{6pc} \label{7.18}
\tau_* = \inf\, \{\, t \ge 0\; \vert\; \varPhi_t \notin (1/\varphi_1^*,
\varphi_1^*)\, \}
\end{equation}
where $\varphi_1^*$ is the unique solution to \eqref{7.16} on
$(1,\infty)$.

To avoid a possible confusion with subscripts we will set $\beta :=
\varphi_1^*$ in the sequel. Thus $\beta \in (0,\infty)$ is the
unique solution to
\begin{equation} \h{7pc} \label{7.19}
\beta - \frac{1}{\beta} + 2\:\! \log \beta = c\:\! \mu^2
\end{equation}
and the stopping time
\begin{equation} \h{6pc} \label{7.20}
\tau = \inf\, \{\, t \ge 0\; \vert\; \varPhi_t \notin (\alpha,
\beta)\, \}
\end{equation}
is optimal in \eqref{7.1} where we set $\alpha = 1/\beta$. These
facts will be used in the subsequent analysis of the optimal
stopping problem \eqref{4.18} when $X$ is three-dimensional.

\section{Properties of the optimal stopping boundaries}
%%%%%%%%%%%%%%%%%%%%%%%%%%%%%%%%%%%%%%%%%%%%%%%%%%%%%%%

In this section we establish the existence of an optimal stopping
time in \eqref{4.18} when the observed process $X$ is
three-dimensional and derive basic properties of the optimal
stopping boundaries. These results will be further refined in
Section 10 below.

\vspace{6pt}

1.\ Looking at \eqref{4.18} we may conclude that the (candidate)
continuation and stopping sets in this problem are respectively
given by
\begin{align} \h{4pc} \label{8.1}
&C = \{\, (\varphi_1,\varphi_2) \in [0,\infty)^2\; \vert\; \hat
V(\varphi_1,\varphi_2) < \hat M(\varphi_1,\varphi_2)\, \} \\[2pt]
\label{8.2} &D = \{\, (\varphi_1,\varphi_2) \in [0,\infty)^2\;
\vert\; \hat V(\varphi_1,\varphi_2) = \hat M(\varphi_1,\varphi_2)
\, \}
\end{align}
where $\hat M$ is defined in \eqref{4.13} above. Recalling
\eqref{5.8} we see that the expectation in \eqref{6.12} defines a
continuous function of the initial point $(\varphi_1,\varphi_2)$ in
$[0,\infty)^2$ for every (bounded) stopping time $\tau$ of $\varPhi$
given and fixed. Taking the infimum over all (bounded) stopping time
$\tau$ of $\varPhi$ we can conclude that the value function $\hat V$
is upper semicontinuous on $[0,\infty)^2$. From \eqref{4.13} and
\eqref{6.11} we see that the loss function $M = \check M \p \hat M$
is continuous and hence lower semicontinuous on $[0,\infty)^2$. It
follows therefore by \cite[Corollary 2.9]{PS} that the first entry
time of the process $\varPhi$ into the closed set $D$ defined by
\begin{equation} \h{6pc} \label{8.3}
\tau_D = \inf\, \{\, t \ge 0\; \vert\; \varPhi_t \in D\, \}
\end{equation}
is optimal in \eqref{6.12}, and hence in \eqref{4.18} as well,
whenever $\PP_{\!\varphi_1,\varphi_2}^0(\tau_D < \infty) = 1$ for
all $(\varphi_1,\varphi_2) \in [0,\infty)^2$. In the sequel we will
establish this and other properties of $\tau_D$ by analysing the
boundary of $D$. We first turn to global properties of the value
function $\hat V$ itself.

\vspace{12pt}

\textbf{Proposition 5.} \emph{For the value function $\hat V$ from
\eqref{4.18} we have}
\begin{align} \h{4pc} \label{8.4}
&(\varphi_1,\varphi_2) \mapsto \hat V(\varphi_1,\varphi_2)\;\;
\text{\emph{is concave on}}\;\; [0,\infty)^2 \\[2pt] \label{8.5}
&(\varphi_1,\varphi_2) \mapsto \hat V(\varphi_1,\varphi_2)\;\;
\text{\emph{is continuous on}}\;\; [0,\infty)^2\, .
\end{align}

\textbf{Proof.} We first show that \eqref{8.4} is satisfied.
Combining \eqref{5.8} with the concavity of the loss function $\hat
M$ from \eqref{4.13} we see that the expectation in \eqref{4.18}
defines a concave function of the initial point
$(\varphi_1,\varphi_2)$ in $[0,\infty)^2$ for every (bounded)
stopping time $\tau$ of $\varPhi$ given and fixed. Taking the
infimum over all (bounded) stopping time $\tau$ of $\varPhi$ we find
that the value function $\hat V$ itself is concave as claimed in
\eqref{8.4} above.

We next show that \eqref{8.5} is satisfied. From the concavity of
$\hat V$ on the open set $(0,\infty)^2$ we can conclude that $\hat
V$ is continuous on $(0,\infty)^2$. Recall that there are concave
functions $F$ defined on a convex subset $S$ of $\R^2$ and taking
values in $\R$, such that the limit of $F(x_n)$ may not exist when
$x_n$ belonging to the interior of $S$ converges to a point $x_0$ at
the boundary of $S$ as $n \rightarrow \infty$. However, if $S$ is
closed then it is well known (and easily verified) that such a
function $F$ must be lower semicontinuous. Applying this implication
to $F = \hat V$ and $S = [0,\infty)^2$ we can conclude that $\hat V$
is lower semicontinuous on $[0,\infty)^2$. At the same time we know
that $\hat V$ is upper semicontinuous (as established following
\eqref{8.2} above) and hence we can conclude that $\hat V$ is
continuous as claimed in \eqref{8.5} above. \hfill $\square$

\vspace{12pt}

2.\ We show that the three straight lines $c_0, c_1, c_2$ defined in
\eqref{6.1}-\eqref{6.3} above are contained in the continuation set
$C$. The proof of this fact uses the Lagrange reformulation
\eqref{6.5} of the optimal stopping problem \eqref{4.18} combined
with the fact that the local times in \eqref{6.5} have a square-root
growth at the three straight lines while the integral in \eqref{6.5}
grows linearly.

\vspace{12pt}

\textbf{Proposition 6.} \emph{The straight lines $c_0, c_1, c_2$
from \eqref{6.1}-\eqref{6.3} are contained in the continuation set
$C$ of the optimal stopping problem \eqref{4.18}.}

\vspace{12pt}

\textbf{Proof.} We claim that
\begin{equation} \h{8pc} \label{8.6}
\EE_{\varphi_1,\varphi_2}^0 \big[ \ell_{t \wedge \tau_{R^c}}^{c_i}
\big] \ge \kappa_i\:\! \sqrt{t}
\end{equation}
for all $t \in (0,t_i)$ with some $\kappa_i > 0$ and $t_i > 0$ for
$i = 0,1,2$ where $\tau_{R^c} = \inf\, \{\, t \ge 0\; \vert\;
(\varPhi_t^1,\varPhi_t^2) \notin R\, \}$ is the first exit time of
$\varPhi$ from a bounded rectangle $R$ containing the given point
$(\varphi_1,\varphi_2) \in c_0 \cup c_1 \cup c_2$ in its interior.
Indeed, this follows by a direct application of Lemma 15 in
\cite{Pe-3} when $(\varphi_1,\varphi_2)$ belongs to $c_0 \cup c_1$,
while the same lemma is applicable to $(\hat \varPhi^1,\hat
\varPhi^2) := (\varPhi^2 \m \varPhi^1,\varPhi^2 \p \varPhi^1)$
obtained by a (bijective) clockwise rotation of
$(\varphi_1,\varphi_2)$ for $45^\circ$ when $(\varphi_1,\varphi_2)$
belongs to $c_2$. Note that the case when $(\varphi_1,\varphi_2) =
(1,1) \in c_0 \cap c_1 \cap c_2$ presents no difficulty as the proof
of Lemma 15 in \cite{Pe-3} extends plainly to cover this case as
well. Having \eqref{8.6} in place we can then proceed as follows.

For $(\varphi_1,\varphi_2) \in c_0 \cup c_1 \cup c_2$ given and
fixed, set $R = [0,2 \varphi_1]\! \times\! [0,2 \varphi_2]$ and
consider the stopping time $\tau := t \wedge \tau_{R^c}$ for $t \in
(0,t_i)$ if $(\varphi_1,\varphi_2)$ belongs to $c_i$ for $i=0,1,2$.
Inserting this $\tau$ under the expectation sign in \eqref{6.5} and
making use of \eqref{8.6} we find that
\begin{equation} \h{4pc} \label{8.7}
\hat V(\varphi_1,\varphi_2) \le (1 \p 2 \varphi_1 \p 2 \varphi_2)\;\!
t - \frac{c}{2}\;\! \kappa_i\;\! \sqrt{t} + \hat M(\varphi_1,\varphi_2)
\end{equation}
for all $t \in (0,t_i)$ if $(\varphi_1,\varphi_2)$ belongs to $c_i$
for $i=0,1,2$. Taking $t$ in \eqref{8.7} sufficiently small we see
that $\hat V(\varphi_1,\varphi_2) < \hat M(\varphi_1,\varphi_2)$
which shows that $(\varphi_1,\varphi_2)$ belongs to $C$ as claimed.
\hfill $\square$

\vspace{6pt}

3.\ The three straight lines $c_0,c_1,c_2$ naturally split the
stopping set $D$ into the three subsets
\begin{align} \h{6pc} \label{8.8}
&D_0 = \{\, (\varphi_1,\varphi_2) \in D\; \vert\; \varphi_1,
\varphi_2 \in [0,1]\, \} \\[2pt] \label{8.9}
&D_1 = \{\, (\varphi_1,\varphi_2) \in D\; \vert\; 1 \le \varphi_1
\le \varphi_2 \, \} \\[2pt] \label{8.10} &D_2 = \{\, (\varphi_1,
\varphi_2) \in D\; \vert\; 1 \le \varphi_2 \le \varphi_1 \, \}\, .
\end{align}
Note that the set $D_0$ is surrounded by the straight lines $c_0$
and $c_1$, the set $D_1$ is surrounded by the straight lines $c_1$
and $c_2$, and the set $D_2$ is surrounded by the straight lines
$c_0$ and $c_2$. Clearly $D = D_0 \cup D_1 \cup D_2$ and the sets
$D_0, D_1, D_2$ are disjoint (see Figure 1 below).

\vspace{12pt}

\textbf{Proposition 7.} \emph{The sets $D_0, D_1, D_2$ are convex.}

\vspace{12pt}

\textbf{Proof.} We will show that the set $D_2$ is convex and the
same arguments can be used to show that the sets $D_0$ and $D_1$ are
convex. For this, let $(\varphi_1',\varphi_2')$ and
$(\varphi_1'',\varphi_2'')$ belonging to $D_2$ and $\lambda \in
(0,1)$ be given and fixed. Firstly, note that
\begin{align} \h{2pc} \label{8.11}
\hat V \big( \lambda (\varphi_1',\varphi_2') \p (1 \m \lambda)
(\varphi_1'',\varphi_2'') \big) &= \hat V \big( \lambda \varphi_1'
\p (1 \m \lambda) \varphi_1'',\lambda \varphi_2' \p (1 \m \lambda)
\varphi_2'' \big) \\[3pt] \notag &\le \hat M \big( \lambda \varphi_1'
\p (1 \m \lambda) \varphi_1'',\lambda \varphi_2' \p (1 \m \lambda)
\varphi_2'' \big) \\[3pt] \notag &= c\:\! \big( 1 \p \lambda
\varphi_2' \p (1 \m \lambda) \varphi_2'' \big)
\end{align}
where we use \eqref{4.13} to infer that $\hat M(\varphi_1,\varphi_2)
= c\:\! (1 \p \varphi_2)$ for $(\varphi_1,\varphi_2)$ belonging to
the subset of $[0,\infty)^2$ surrounded by $c_0$ and $c_2$.
Secondly, using that $\hat V$ is concave on $[0,\infty)^2$ as
established in \eqref{8.4} above, we find that
\begin{align} \h{2pc} \label{8.12}
\hat V \big( \lambda (\varphi_1',\varphi_2') \p (1 \m \lambda)
(\varphi_1'',\varphi_2'') \big) &\ge \lambda \hat V(\varphi_1',
\varphi_2') \p (1 \m \lambda) \hat V(\varphi_1'',\varphi_2'')
\\[3pt] \notag &=\lambda\:\! \hat M(\varphi_1',\varphi_2') \p
(1  \m \lambda) \hat M(\varphi_1'',\varphi_2'') \\[3pt] \notag
&= c\:\! \big( 1 \p \lambda \varphi_2' \p (1 \m \lambda)
\varphi_2'' \big)
\end{align}
where in the first equality we use that $(\varphi_1',\varphi_2')$
and $(\varphi_1'',\varphi_2'')$ belong to $D_2 \subseteq D$.
Combining \eqref{8.11} and \eqref{8.12} we see that $\hat V \big(
\lambda (\varphi_1',\varphi_2') \p (1 \m \lambda)
(\varphi_1'',\varphi_2'') \big) = \hat M \big( \lambda
(\varphi_1',\varphi_2') \p (1 \m \lambda) (\varphi_1'',\varphi_2'')
\big)$ showing that $\lambda (\varphi_1',\varphi_2') \p (1 \m
\lambda) (\varphi_1'',\varphi_2'')$ belongs to $D_2$ as needed.
\hfill $\square$

\vspace{12pt}

4.\ To describe the shape of the stopping sets $D_0, D_1, D_2$ we
may recall from Section 7 that the subsets $\big( [0,1/\beta] \cup
[\beta,\infty) \big)\! \times\! \{0\}$ and $\{0\}\! \times\! \big(
[0,1/\beta] \cup [\beta,\infty) \big)$ of $[0,\infty)^2$ are
contained in $D$ where $\beta \in (1,\infty)$ solves \eqref{7.19}
uniquely. Symmetry arguments to be addressed shortly below show that
it is sufficient to focus on the set $D_2$ as the conclusions will
directly extend to the sets $D_0$ and $D_1$ as well. Moving from the
straight lines $c_0$ and $c_2$ in $C$ to the right, let us formally
define the (least) boundary between $C$ and $D_2$ by setting
\begin{equation} \h{5pc} \label{8.13}
b_2(\varphi_2) = \inf\, \{\, \varphi_1 > 1 \vee \varphi_2\; \vert\;
(\varphi_1,\varphi_2) \in D_2\, \}
\end{equation}
for $\varphi_2 \in [0,\infty)$. Clearly the infimum in \eqref{8.13}
is attained since $D_2$ is closed. We now show that $b_2$
constitutes the entire boundary of $D_2$ in $[0,\infty)^2$ (see
Figure 1 below).

\vspace{12pt}

\textbf{Proposition 8.} \emph{The mapping $\varphi_2 \mapsto
b_2(\varphi_2)$ is finite valued on $[0,\infty)$ and we have
\begin{equation} \h{5pc} \label{8.14}
D_2 = \{\, (\varphi_1,\varphi_2) \in [0,\infty)^2\; \vert\;
\varphi_1 \ge b_2(\varphi_2)\, \}
\end{equation}
with $b_2(0) = \beta \in (1,\infty)$ and $b_2(\varphi_2) \rightarrow
\infty$ as $\varphi_2 \rightarrow \infty$.}

\vspace{12pt}

\textbf{Proof.} To derive \eqref{8.14} we show that
\begin{equation} \h{6pc} \label{8.15}
(\varphi_1',\varphi_2') \in D_2 \Rightarrow (\varphi_1'',
\varphi_2') \in D_2
\end{equation}
for all $\varphi_1'' \ge \varphi_1'$. For this, recall from
\eqref{8.4} that $\varphi_1 \mapsto \hat V(\varphi_1,\varphi_2')$ is
concave on $[0,\infty)$ while $\varphi_1 \mapsto \hat
M(\varphi_1,\varphi_2') = c\:\! (1 \p \varphi_2')$ is constant for
$\varphi_1 \ge \varphi_2' \ge 1$. Hence if $\hat
V(\varphi_1',\varphi_2') = \hat M (\varphi_1',\varphi_2')$ due to
$(\varphi_1',\varphi_2') \in D_2$ with $\hat
V(\varphi_1'',\varphi_2') < \hat M (\varphi_1'',\varphi_2')$ meaning
that $(\varphi_1'',\varphi_2') \notin D_2$ for some $\varphi_1'' >
\varphi_1'$, then $\hat V(\varphi_1,\varphi_2')$ must converge to
$-\infty$ as $\varphi_1$ converges to $\infty$. This however
contradicts the fact that $\hat V$ is non-negative and hence
\eqref{8.15} must hold as claimed. Combining \eqref{8.13} and
\eqref{8.15} we see that \eqref{8.14} is satisfied as claimed.

To establish that $b_2$ is finite valued we first show that
\begin{equation} \h{6pc} \label{8.16}
(\varphi_1',\varphi_2') \in D_2 \Rightarrow (\varphi_1',
\varphi_2'') \in D_2
\end{equation}
for all $\varphi_2'' \in [0,\varphi_2']$ when $\varphi_1' \ge
\beta$. (Note that the latter inequality cannot be omitted and
\eqref{8.16} may fail when $\varphi_1' < \beta$ as we will see in
Section 10 below.) For this, recall from \eqref{8.4} that $\varphi_2
\mapsto \hat V(\varphi_1',\varphi_2)$ is concave on $[0,\infty)$
while $\varphi_2 \mapsto \hat M(\varphi_1',\varphi_2) = c\:\! (1 \p
\varphi_2)$ is linear for $\varphi_2 \in [1,\varphi_1']$. By the
results of Section 7 we know that $(\varphi_1',0)$ belongs to $D_2$
so that $\hat V(\varphi_1',0) = \hat M(\varphi_1',0)$ when
$\varphi_1' \ge \beta$. Hence if $\hat V(\varphi_1',\varphi_2') =
\hat M(\varphi_1',\varphi_2')$ due to $(\varphi_1',\varphi_2') \in
D_2$ then $\hat V(\varphi_1',\varphi_2) = \hat
M(\varphi_1',\varphi_2)$ for all $\varphi_2 \in [0,\varphi_2']$.
This shows that \eqref{8.16} holds as claimed.

From \eqref{8.16} we see that if $b_2(\varphi_2') \ge \beta$ for
some $\varphi_2' > 0$ then $\varphi_2 \mapsto b_2(\varphi_2)$ is
increasing on $[\varphi_2',\infty)$. In particular, this means that
if $b_2(\varphi_2'') = \infty$ for some $\varphi_2'' > 0$ then
$b_2(\varphi_2) = \infty$ for all $\varphi_2 \ge \varphi_2''$. We
will now use this fact to show that $b_2$ is finite valued as
claimed.

Assuming that $b_2(\varphi_2'') = \infty$ for some $\varphi_2''
> 0$, and fixing $b > a > \varphi_2''$, it follows from the previous
argument that the rectangle $R_N = (N,\infty)\! \times\! (a,b)$ is
contained in $C$ for every $N \ge N_0$ with some $N_0 \ge 1$ large
enough. Consider the stopping time
\begin{equation} \h{6pc} \label{8.17}
\tau_{R_N^c}^{\varphi_1,\varphi_2} = \inf\, \{\, t \ge 0\; \vert\;
(\varPhi_t^{\varphi_1},\varPhi_t^{\varphi_2}) \notin R_N\, \}
\end{equation}
for $(\varphi_1,\varphi_2) \in R_N$. Since $R_N \subseteq C$ we see
that $\tau_{R_N^c}^{\varphi_1,\varphi_2} \le
\tau_D^{\varphi_1,\varphi_2}$ and hence it follows that
\begin{align} \h{1pc} \label{8.18}
c\:\! (1 \p \varphi_2) &\ge \hat M(\varphi_1,\varphi_2) \ge
\hat V(\varphi_1,\varphi_2) = \EE_{\varphi_1,\varphi_2}^0 \Big[
\int_0^{\tau_D}\! \big( 1 \p \varPhi_t^1 \p \varPhi_t^2 \big)\,
dt + \hat M \big( \varPhi_{\tau_D}^1, \varPhi_{\tau_D}^2 \big)
\:\! \Big] \\[-2pt] \notag &\ge \EE_0 \Big[ \int_0^{\tau_{R_N^c}^
{\varphi_1,\varphi_2}}\!\! \varPhi_t^{\varphi_1}\, dt\, \Big]
\ge N\, \EE_0 \big[ \tau_{R_N^c}^{\varphi_1,\varphi_2} \big]
\end{align}
for all $N \ge N_0$. Noting that $\EE_0 \big[
\tau_{R_N^c}^{\varphi_1,\varphi_2} \big] \rightarrow \EE_0 \big[
\tau_{(a,b)^c}^{\varphi_2} \big]$ where $\tau_{(a,b)^c}^{\varphi_2}
= \inf\, \{\, t \ge 0\; \vert\; \varPhi_t^{\varphi_2} \notin (a,b)\,
\}$ as $N \rightarrow \infty$, we see from \eqref{8.18} that
\begin{equation} \h{8pc} \label{8.19}
c\:\! (1 \p \varphi_2) \ge \frac{N}{2}\, \EE_0 \big[ \tau_{(a,b)^c}^
{\varphi_2} \big]
\end{equation}
for all $N \ge N_1$ with some $N_1 \ge 1$ large enough. Letting $N
\rightarrow \infty$ and using that $\EE_0 \big[
\tau_{(a,b)^c}^{\varphi_2} \big]
> 0$ we obtain a contradiction. Thus there is no $\varphi_2'' > 0$
such that $b_2(\varphi_2'') = \infty$ and hence $b_2$ is finite
valued as claimed.

Finally, the fact that $b_2(0) = \beta \in (1,\infty)$ was
established in Section 7 above. Moreover, since $b_2(\varphi_2) >
\varphi_2$ for all $\varphi_2 > 0$ due to $c_0$ and $c_2$ being
contained in $C$, we see that $b_2(\varphi_2) \rightarrow \infty$ as
$\varphi_2 \rightarrow \infty$ and the proof is complete. \hfill
$\square$

\vspace{12pt}

\textbf{Proposition 9.} \emph{The mapping $\varphi_2 \mapsto
b_2(\varphi_2)$ is convex and continuous on $[0,\infty)$.}

\vspace{12pt}

\textbf{Proof.} Convexity of the mapping $\varphi_2 \mapsto
b_2(\varphi_2)$ on $[0,\infty)$ follows from the convexity of the
stopping set $D_2$ as established in Proposition 7 above. Hence the
mapping $\varphi_2 \mapsto b_2(\varphi_2)$ is continuous on
$(0,\infty)$ while $b_2$ cannot make a jump at $0$ due to the fact
that the stopping set $D_2$ is closed. This completes the proof.
\hfill $\square$

\vspace{12pt}

We will show in Section 10 below that $b_2(\varphi_2) < \beta$ for
$\varphi_2 \in (0,\kappa)$ with $\kappa > 0$ such that $b_2(\kappa)
= \beta$. This fact combined with the convexity of $b_2$ on
$[0,\infty)$ means that the mapping $\varphi_2 \mapsto
b_2(\varphi_2)$ is (firstly) decreasing on $[0,\kappa']$ and (then)
increasing on $[\kappa',\infty)$ with some $\kappa' \in (0,\kappa)$.
In addition to these facts about $b_2$ around zero we will conclude
this section by evaluating the asymptotic behaviour of $b_2$ at
infinity. Before we do that we will turn to the remaining two
stopping sets $D_0$ and $D_1$ including their boundaries.

\vspace{6pt}

5.\ Symmetry arguments enable us to extend the setting and results
of Proposition 8 and Proposition 9 from the stopping set $D_2$ to
the remaining two stopping sets $D_0$ and $D_1$. For this, recall
from \eqref{4.3} that $\varPhi^1 = \varPi^1/\varPi^0$ and $\varPhi^2
= \varPi^2/\varPi^0$. Since $\varPi^0, \varPi^1, \varPi^2$ play a
symmetric role in the optimal stopping problem \eqref{3.8} we see
that any permutation of the three coordinates should yield the same
result. There are two generic permutations which generate all the
others (six in total). The first generic permutation is obtained by
swapping $\varPi^1$ and $\varPi^2$ while keeping $\varPi^0$ intact.
This yields $\varPhi^1 = \varPi^1/\varPi^0 \sim \varPi^2/\varPi^0 =
\varPhi^2$ and $\varPhi^2 = \varPi^2/\varPi^0 \sim \varPi^1/\varPi^0
= \varPhi^1$ showing that
\begin{equation} \h{5pc} \label{8.20}
(\varphi_1,\varphi_2) \in \partial C \Longleftrightarrow (\varphi_2,
\varphi_1) \in \partial C
\end{equation}
where $\partial C$ can also be replaced by $C$ or $D$. The second
generic permutation is obtained by swapping $\varPi^0$ and
$\varPi^1$ while keeping $\varPi^2$ intact. This yields $\varPhi_1 =
\varPi^1/\varPi^0 \sim \varPi^0/\varPi^1 = 1/\varPhi^1$ and
$\varPhi_2 = \varPi^2/\varPi^0 \sim \varPi^2/\varPi^1 =
\varPhi^2/\varPhi^1$ showing that
\begin{equation} \h{5pc} \label{8.21}
(\varphi_1,\varphi_2) \in \partial C \Longleftrightarrow \Big(
\frac{1}{\varphi_1}, \frac{\varphi_2}{\varphi_1} \Big) \in \partial C
\end{equation}
where $\partial C$ can also be replaced by $C$ or $D$. The remaining
four equivalencies can be obtained by combining \eqref{8.20} and
\eqref{8.21}. For example, applying first \eqref{8.20} and then
\eqref{8.21} we find that $(\varphi_1,\varphi_2) \in \partial C
\Leftrightarrow (1/\varphi_2, \varphi_1/\varphi_2) \in \partial C$
(where $\partial C$ can also be replaced by $C$ or $D$ as above)
which is obtained by swapping $\varPi^0$ and $\varPi^2$ while
keeping $\varPi^1$ intact.

\vspace{6pt}

6.\ Having understood the symmetry relations we now move to
extending the setting and results of Proposition 8 and Proposition 9
from $D_2$ to $D_0$ and $D_1$. We first address the case of $D_1$
which in view of \eqref{8.20} is a mirror image of $D_2$ across the
main diagonal in $[0,\infty)^2$. In analogy with \eqref{8.13} we
thus define the (least) boundary between $C$ and $D_1$ by setting
\begin{equation} \h{5pc} \label{8.22}
b_1(\varphi_1) = \inf\, \{\, \varphi_2 > 1 \vee \varphi_1\; \vert\;
(\varphi_1,\varphi_2) \in D_1\, \}
\end{equation}
for $\varphi_1 \in [0,\infty)$. Clearly the infimum in \eqref{8.22}
is attained since $D_1$ is closed. Similarly to $b_2$ and $D_2$
above we now show that $b_1$ constitutes the entire boundary of
$D_1$ in $[0,\infty)^2$ (see Figure 1 below).

\vspace{12pt}

\textbf{Proposition 10.} \emph{The mapping $\varphi_1 \mapsto
b_1(\varphi_1)$ is finite valued on $[0,\infty)$ and we have
\begin{equation} \h{5pc} \label{8.23}
D_1 = \{\, (\varphi_1,\varphi_2) \in [0,\infty)^2\; \vert\;
\varphi_2 \ge b_1(\varphi_1)\, \}
\end{equation}
with $b_1(0) = \beta \in (1,\infty)$ and $b_1(\varphi_1) \rightarrow
\infty$ as $\varphi_1 \rightarrow \infty$.}

\vspace{12pt}

\textbf{Proof.} This can be derived in exactly the same way as in
Proposition 8 above. Alternatively Proposition 10 also follows
directly from Proposition 8 using the symmetry relation \eqref{8.20}
which shows that $D_1$ is a mirror image of $D_2$ across the main
diagonal in $[0,\infty)^2$. \hfill $\square$

\vspace{12pt}

\textbf{Proposition 11.} \emph{The mapping $\varphi_1 \mapsto
b_1(\varphi_1)$ is convex and continuous on $[0,\infty)$.}

\vspace{12pt}

\textbf{Proof.} This can be derived in exactly the same way as in
Proposition 9 above. Alternatively Proposition 11 also follows
directly from Proposition 9 using the symmetry relation \eqref{8.20}
which shows that $b_1$ coincides with $b_2$ on $[0,\infty)$. \hfill
$\square$

\vspace{12pt}

Despite the fact that the functional rules of $b_1$ and $b_2$
coincide on $[0,\infty]$, we will still keep their different
subscripts $1$ and $2$ in place to account for different arguments
in $(\varphi_1,b_1(\varphi_1))$ and $(b_2(\varphi_2),\varphi_2)$ for
$\varphi_1 \ge 0$ and $\varphi_2 \ge 0$ respectively.

\vspace{12pt}

6.\ We next address the case of $D_0$ which in view of \eqref{8.21}
can similarly be linked to the case of $D_2$ in a one-to-one way.
Moving from the point $(1,1) \in C$ down to the point $(0,0) \in
D_0$ along the main diagonal in $[0,1]^2$, we know that there exists
the (first) point $(\gamma,\gamma)$ that belongs to $D_0$.
Equivalently $\gamma$ can also be formally defined by
\begin{equation} \h{5pc} \label{8.24}
\gamma = \sup\, \{\, \varphi \in [0,1]\; \vert\; (\varphi,\varphi)
\in D_0\, \}\, .
\end{equation}
Clearly the supremum in \eqref{8.24} is attained since $D_0$ is
closed and we have $\gamma \in (0,1)$ since $(1,1) \in C$. Similarly
to \eqref{8.13} and \eqref{8.22} we then define the (least) upper
boundary between $C$ and $D_0$ by setting
\begin{equation} \h{5pc} \label{8.25}
b_0^1(\varphi_1) = \sup\, \{\, \varphi_2 \in [0,1]\; \vert\;
(\varphi_1,\varphi_2) \in D_0\, \}
\end{equation}
for $\varphi_1 \in [0,\gamma]$, and the (least) lower boundary
between $C$ and $D_0$ by setting
\begin{equation} \h{5pc} \label{8.26}
b_0^2(\varphi_2) = \sup\, \{\, \varphi_1 \in [0,1]\; \vert\;
(\varphi_1,\varphi_2) \in D_0\, \}
\end{equation}
for $\varphi_2 \in [0,\gamma]$. Clearly the suprema in \eqref{8.25}
and \eqref{8.26} are attained since $D_0$ is closed. In view of
\eqref{8.20}, it is clear that the graphs of $b_0^1$ and $b_0^2$ are
mirror images of each other across the main diagonal in
$[0,\gamma]^2$, so that $b_0^1 = b_0^2$ on $[0,\gamma]$ and we set
\begin{equation} \h{7pc} \label{8.27}
b_0(\varphi) := b_0^1(\varphi) = b_0^2(\varphi)
\end{equation}
for $\varphi \in [0,\gamma]$. Similarly to Proposition 8 and
Proposition 10 above, we now show that $b_0$ can be used to describe
the entire boundary of $D_0$ in $[0,\infty)^2$\! (see Figure 1
below).

\vspace{12pt}

\textbf{Proposition 12.} \emph{The following identity holds
\begin{equation} \h{1.5pc} \label{8.28}
D_0 = \{\, (\varphi_1,\varphi_2) \in [0,\gamma]^2\; \vert\; \varphi_1
\le \varphi_2 \le b_0(\varphi_1)\;\; \text{or}\;\; \varphi_2 \le
\varphi_1 \le b_0(\varphi_2)\, \}
\end{equation}
with $b_0(0) = 1/\beta$ and $b_0(\gamma) = \gamma \in (1/\beta,1)$.
The mapping $\varphi \mapsto b_0(\varphi)$ is concave and conti-
nuous on $[0,\gamma]$.}

\vspace{6pt}

\textbf{Proof.} All claims follow by convexity (and closeness) of
$D_0$ established in Proposition 7 above combined with the symmetry
relations \eqref{8.20} and \eqref{8.21}. The latter symmetry
relation links $D_0$ to $D_2$ in a one-to-one way and this enables
us to conclude that $b_0(0) = 1/\beta$ as claimed. The final claim
$b_0(\gamma) = \gamma \in (1/\beta,1)$ is evident from
\eqref{8.24}-\eqref{8.26} above. \hfill $\square$

\vspace{12pt}

The one-to-one correspondence between $D_0$ and $D_2$ obtained by
the symmetry relation \eqref{8.21} enables us to transfer the facts
stated following Proposition 9 above from $D_2$ to $D_0$. In
particular, this yields that the mapping $\varphi \mapsto
b_0(\varphi)$ is (firstly) increasing on $[0,\delta]$ and (then)
decreasing on $[\delta,\gamma]$ for some $\delta \in (0,\gamma)$
(see Figure 1 below).

\vspace{12pt}

7.\ Another consequence of the one-to-one correspondence between
$D_0$ and $D_2$ (and hence $D_1$ as well) is the possibility to
describe the asymptotic behaviour of $b_1$ and $b_2$ at infinity.

\vspace{12pt}

\textbf{Proposition 13.} \emph{We have
\begin{equation} \h{6pc} \label{8.29}
\lim_{\varphi_1 \rightarrow \infty} \frac{b_1(\varphi_1)}{\varphi_1}
= \lim_{\varphi_2 \rightarrow \infty} \frac{b_2(\varphi_2)}{\varphi_2}
= \beta\, .
\end{equation}}

\textbf{Proof.} The first equality follows by the symmetry relation
\eqref{8.20} implying that $b_1$ coincides with $b_2$ on
$[0,\infty)$ so that it is enough to establish the second equality
in \eqref{8.29}. For this, note that the symmetry relation
\eqref{8.21} yields
\begin{equation} \h{6pc} \label{8.30}
(\varphi_1,b_0(\varphi_1)) \in \partial C \Longleftrightarrow
\Big( \frac{1}{\varphi_1}, \frac{b_0(\varphi_1)}{\varphi_1}
\Big) \in \partial C
\end{equation}
for $\varphi_1 \in [0,\gamma]$. Note also that
$(\varphi_1,b_0(\varphi_1)) \in \partial D_0$ tends to $(0,1/\beta)$
and $(1/\varphi_1,b_0(\varphi_1)/\varphi_1) \in
\partial D_2$ tends to $(\infty,\infty)$ as $\varphi_1 \rightarrow
0$. The fact that the point $(1/\varphi_1,b_0(\varphi_1)/\varphi_1)$
belongs to $\partial C \cap \partial D_2$ means that this point can
be identified with $(b_2(\varphi_2),\varphi_2)$ for some $\varphi_2
> 0$ with $\varphi_2 \rightarrow \infty$ as $\varphi_1 \rightarrow
0$. This shows that
\begin{equation} \h{6pc} \label{8.31}
\frac{b_2(\varphi_2)}{\varphi_2} = \frac{\frac{1}{\varphi_1}}
{\frac{b_0(\varphi_1)}{\varphi_1}} = \frac{1}{b_0(\varphi_1)}
\rightarrow \frac{1}{\frac{1}{\beta}} = \beta
\end{equation}
as $\varphi_2 \rightarrow \infty$. This establishes \eqref{8.29} and
the proof is complete. \hfill $\square$

\vspace{12pt}

We will continue our study of the sets $D_0, D_1, D_2$ in Section 10
below.

\section{Smooth fit}
%%%%%%%%%%%%%%%%%%%%

In this section we show that the value function $\hat V$ from
\eqref{4.18} satisfies the smooth fit condition at the optimal
stopping boundaries $b_0, b_1, b_2$. A key point in the proof is
based upon the fact that the boundary points are \emph{Green
regular} for $D_0, D_1, D_2$ in the sense that the first entry time
$\tau_{D_i}^{\varphi_1^n,\varphi_2^n}$ of
$(\varPhi^{\varphi_1^n},\varPhi^{\varphi_2^n})$ into $D_i$ satisfies
\begin{equation} \h{9pc} \label{9.1}
\tau_{D_i}^{\varphi_1^n,\varphi_2^n} \rightarrow 0
\end{equation}
with $\PP_{\!0}$\!-probability one whenever
$(\varphi_1^n,\varphi_2^n)$ from $C$ tends to
$(\varphi_1,\varphi_2)$ at the boundary $\partial C \cap D_i$ for
$i=0,1,2$ as $n \rightarrow \infty$. The Green regularity follows
from the fact that the boundary points are \emph{probabilistically
regular} for $D_0, D_1, D_2$ in the sense that
$\PP_{\!\varphi_1,\varphi_2}^0( \tau_{D_i}\! =\! 0) = 1$ for every
$(\varphi_1,\varphi_2)$ at the boundary $\partial C \cap D_i$ for
$i=0,1,2$ combined with the fact that the process
$(\varPhi^1,\varPhi^2)$ is strong Feller which is evident from
\eqref{5.8} above (cf.\ \cite[Section 3]{DePe}). The probabilistic
regularity is a consequence of the fact that the sets $D_i$ are
convex (as established in Proposition 7 above) so that in view of
\eqref{5.8} each boundary point from $\partial C \cap D_i$ satisfies
Zaremba's cone condition for $D_i$ with $i=0,1,2$ (see e.g.\
\cite[Theorem 3.2, p.\ 250]{KS}). These facts establish \eqref{9.1}
and we can now state the main result of this section.

\vspace{12pt}

\textbf{Proposition 14 (Smooth fit).} \emph{For the value function $\hat
V$ from \eqref{4.18} we have
\begin{align} \h{7pc} \label{9.2}
&\hat V_{\varphi_1}(\varphi_1,\varphi_2) = \hat M_{\varphi_1}(\varphi_1,
\varphi_2) \\[2pt] \label{9.3}&\hat V_{\varphi_2}(\varphi_1,\varphi_2)
= \hat M_{\varphi_2}(\varphi_1, \varphi_2)
\end{align}
for all $(\varphi_1,\varphi_2) \in \partial C \cap D_i$ with
$i=0,1,2$.}

\vspace{12pt}

\textbf{Proof.} We will establish \eqref{9.2} and \eqref{9.3} for
$D_2$ and similar arguments can be used for $D_0$ and $D_1$. For
this, let $\varphi_1 = b_2(\varphi_2)$ with $\varphi_2>0$ be given
and fixed in the sequel.

\vspace{6pt}

1.\ We show that \eqref{9.2} holds. For this, we first note that
\begin{equation} \h{1pc} \label{9.4}
\liminf_{h \downarrow 0}\, \frac{\hat V(\varphi_1 \m h,\varphi_2) \m
\hat V(\varphi_1,\varphi_2)}{-h} \ge \liminf_{h \downarrow 0}\,
\frac{\hat M(\varphi_1 \m h,\varphi_2) \m \hat M(\varphi_1,
\varphi_2)}{-h} = 0
\end{equation}
since $\hat V(\varphi_1 \m h,\varphi_2) \le \hat M(\varphi_1 \m
h,\varphi_2)$ and $\hat V(\varphi_1,\varphi_2) = \hat
M(\varphi_1,\varphi_2)$ with $\varphi_1' \mapsto \hat
M(\varphi_1',\varphi) = c\:\! (1 +$ $\varphi_2)$ being constant for
$\varphi_1' > \varphi_2 \ge 1$. We next show that
\begin{equation} \h{5pc} \label{9.5}
\limsup_{h \downarrow 0}\, \frac{\hat V(\varphi_1 \m h,\varphi_2) \m
\hat V(\varphi_1,\varphi_2)}{-h} \le 0\, .
\end{equation}
For this, let $\tau_D^{\varphi_1-h,\varphi_2}$ denote the first
entry time of $(\varPhi^{\varphi_1-h},\varPhi^{\varphi_2})$ into $D$
for $h>0$ given and fixed. Since $\tau_D^{\varphi_1-h,\varphi_2}$ is
optimal for $\hat V(\varphi_1 \m h,\varphi_2)$ we find by
\eqref{5.8} that
\begin{align} \h{-0.5pc} \label{9.6}
\hat V(\varphi_1 \m h,\varphi_2) \m \hat V(\varphi_1,\varphi_2) &\ge
\EE_0 \Big[ \int_0^{\tau_D^{\varphi_1-h,\varphi_2}}\!\!\!\!\!\!\!\!
(1 \p \varPhi_t^{\varphi_1-h}\! \p \varPhi_t^{\varphi_2})\, dt +
\hat M \big( \varPhi_{\!\tau_D^{\varphi_1-h,\varphi_2}}^{\varphi_1-h},
\varPhi_{\!\tau_D^{\varphi_1-h,\varphi_2}}^{\varphi_2} \big) \Big]
\\[-3pt] \notag &\h{13pt}- \EE_0 \Big[ \int_0^{\tau_D^{\varphi_1-h,
\varphi_2}}\!\!\!\!\!\! (1 \p \varPhi_t^{\varphi_1}\! \p \varPhi_t^
{\varphi_2})\, dt + \hat M \big( \varPhi_{\!\tau_D^{\varphi_1-h,
\varphi_2}}^{\varphi_1},\varPhi_{\!\tau_D^{\varphi_1-h,\varphi_2}}^
{\varphi_2} \big) \Big] \\[-1pt] \notag &= \EE_0 \Big[ \int_0^{\tau_
D^{\varphi_1-h,\varphi_2}}\!\!\!\!\!\!\!\! (-h\:\! \varPhi_t^1)\, dt
- c\:\! h\:\! \varPhi_{\!\tau_D^{\varphi_1-h,\varphi_2}}^1\, I \big(
\tau_D^{\varphi_1-h,\varphi_2} \ne \tau_{D_2}^{\varphi_1-h,\varphi_2}
\big) \Big]
\end{align}
for all $h \in (0,h_0)$ with some $h_0>0$ sufficiently small, where
in the final equality we use \eqref{6.4} combined with the two
implications \eqref{9.7} below which we motivate and derive first.

Recalling \eqref{6.4} and definitions of $\Delta_0, \Delta_1,
\Delta_2$ stated afterwards, we claim that
\begin{align} \h{2pc} \label{9.7}
\big( \varPhi_{\!\tau_D^{\varphi_1-h,\varphi_2}}^{\varphi_1-h},
\varPhi_{\!\tau_D^{\varphi_1-h,\varphi_2}}^{\varphi_2} \big) \in
\partial C \cap D_i \Longrightarrow \big( \varPhi_{\!\tau_D^{
\varphi_1-h,\varphi_2}}^{\varphi_1},\varPhi_{\!\tau_D^{\varphi_
1-h,\varphi_2}}^{\varphi_2} \big) \in \Delta_i
\end{align}
for $h \in (0,h_0)$ with some $h_0>0$ sufficiently small and
$i=0,1$.

To show \eqref{9.7} for $i=0$ recall that $c_0$ and $c_1$ are
contained in $C$ so that the continuous curve $b_0$ stays away from
the straight line $c_0$ in particular. Setting $\tau_h :=
\tau_{D_0}^{\varphi_1-h,\varphi_2}$ to simplify the notation
throughout this shows that there exists $\delta>0$ sufficiently
small such that the right-hand side in \eqref{9.7} with $i=0$
implies that $(\varphi_1 \m h)\;\! \varPhi_{\tau_h}^1 \le 1 \m
\delta$ for $h \in (0,\varphi_1)$. This imp- lies that
$\varPhi_{\tau_h}^1 \le (1 \m \delta)/(\varphi_1 \m h_0)$ for all $h
\in (0,h_0)$ with $h_0 \in (0,\varphi_1)$ given and fixed. It
follows that $\varPhi_{\tau_h}^{\varphi_1} = \varphi_1\:\!
\varPhi_{\tau_h}^1 = (\varphi_1 \m h)\:\! \varPhi_{\tau_h}^1 + h\:\!
\varPhi_{\tau_h}^1 \le 1 \m \delta + h_0\:\! (1 \m
\delta)/(\varphi_1 \m h_0) = [(1 \m \delta)(1 \p h_0)/(\varphi_1 \m
h_0)] \le 1$ if we choose $h_0>0$ small enough. This shows that
\eqref{9.7} holds for $i=0$ as claimed.

To show \eqref{9.7} for $i=1$ set $\tau_h :=
\tau_{D_1}^{\varphi_1-h,\varphi_2}$ to simplify the notation
throughout and note that \eqref{8.29} shows that there exists
$\delta>0$ sufficiently small such that the right-hand side in
\eqref{9.7} with $i=1$ implies that $(\varphi_1 \m h)\;\!
\varPhi_{\tau_h}^1 \le (1 \m \delta)\:\! \varphi_2\;\!
\varPhi_{\tau_h}^2$ for $h \in (0,\varphi_1)$. This implies that
$\varPhi_{\tau_h}^1 \le [(1 \m \delta)/(\varphi_1 \m h_0)]\:\!
\varphi_2\;\! \varPhi_{\tau_h}^2$ for all $h \in (0,h_0)$ with $h_0
\in (0,\varphi_1)$ given and fixed. It follows that
$\varPhi_{\tau_h}^{\varphi_1} = \varphi_1\:\! \varPhi_{\tau_h}^1 =
(\varphi_1 \m h)\:\! \varPhi_{\tau_h}^1 \p h\;\! \varPhi_{\tau_h}^1
\le (1 \m \delta)\;\! \varphi_2\;\! \varPhi_{\tau_h}^2 + h\:\! [(1
\m \delta)/(\varphi_1 \m h_0)]\;\! \varphi_2\;\! \varPhi_{\tau_h}^2$
$\le [(1 \m \delta)(1 \p h_0)/(\varphi_1 \m h_0)]
\varPhi_{\tau_h}^{\varphi_2} \le \varPhi_{\tau_h}^{\varphi_2}$ if we
choose $h_0>0$ small enough. This shows that \eqref{9.7} holds for
$i=1$ as claimed.

Making now use of \eqref{6.4} and \eqref{9.7} in the middle term of
\eqref{9.6} above, upon noting that $\tau_D^{\varphi_1-h,\varphi_2}$
always equals one among $\tau_{D_0}^{\varphi_1-h,\varphi_2}$,
$\tau_{D_1}^{\varphi_1-h,\varphi_2}$,
$\tau_{D_2}^{\varphi_1-h,\varphi_2}$ respectively, we see that the
final equality in \eqref{9.6} holds as claimed. Dividing both sides
of \eqref{9.6} by $-h$ we obtain
\begin{equation} \h{1pc} \label{9.8}
\frac{\hat V(\varphi_1 \m h,\varphi_2) \m \hat V(\varphi_1,\varphi_2)}
{-h} \le \EE_0 \Big[ \int_0^{\tau_D^{\varphi_1-h,\varphi_2}}\!\!\!\!\!
\!\!\! \varPhi_t^1\, dt + c\;\! \varPhi_{\!\tau_D^{\varphi_1-h,\varphi_
2}}^1\, I \big( \tau_D^{\varphi_1-h,\varphi_2} \ne \tau_{D_2}^{\varphi_
1-h,\varphi_2} \big) \Big]
\end{equation}
for all $h \in (0,h_0)$. Letting $h \downarrow 0$ and using that the
right-hand side in \eqref{9.8} tends to zero by \eqref{9.1} and the
continuity of $\hat V$ we see that \eqref{9.5} holds as claimed.
Combining \eqref{9.4} and \eqref{9.5} with the fact that $\hat
M_{\varphi_1}(\varphi_1,\varphi_2)=0$ we see that \eqref{9.2} holds
as claimed.

\vspace{6pt}

2.\ We show that \eqref{9.3} holds. For this, we first note that
\begin{align} \h{1pc} \label{9.9}
&\liminf_{h \downarrow 0}\, \frac{\hat V(\varphi_1,\varphi_2 \m h)
\m \hat V(\varphi_1,\varphi_2)}{-h} \ge \liminf_{h \downarrow 0}
\, \frac{\hat M(\varphi_1,\varphi_2 \m h) \m \hat M(\varphi_1,
\varphi_2)}{-h} = c \\ \label{9.10} &\limsup_{h \downarrow 0}\,
\frac{\hat V(\varphi_1,\varphi_2 \p h) \m \hat V(\varphi_1,
\varphi_2)}{h} \le \limsup_{h \downarrow 0}\, \frac{\hat M(
\varphi_1,\varphi_2 \p h) \m \hat M(\varphi_1,\varphi_2)}{h}
= c
\end{align}
depending on whether $(\varphi_1,\varphi_2 \m h)$ or
$(\varphi_1,\varphi_2 \p h)$ belongs to $C$ for $h>0$ respectively.
In \eqref{9.8} and \eqref{9.9} we use that $\hat
V(\varphi_1,\varphi_2 \mp h) \le \hat M(\varphi_1,\varphi_2 \mp h)$
and $\hat V(\varphi_1,\varphi_2) = \hat M(\varphi_1,\varphi_2)$ with
$\varphi_2' \mapsto \hat M(\varphi_1,\varphi_2') = c\:\! (1 \p
\varphi_2')$ being linear for $1 \le \varphi_2' < \varphi_1$. We
next show that
\begin{align} \h{6pc} \label{9.11}
&\limsup_{h \downarrow 0}\, \frac{\hat V(\varphi_1,\varphi_2 \m h)
\m \hat V(\varphi_1,\varphi_2)}{-h} \le c \\ \label{9.12} &\liminf_
{h \downarrow 0}\, \frac{\hat V(\varphi_1,\varphi_2 \p h)
\m \hat V(\varphi_1,\varphi_2)}{h} \ge c
\end{align}
depending on whether $(\varphi_1,\varphi_2 \m h)$ or
$(\varphi_1,\varphi_2 \p h)$ belongs to $C$ for $h>0$ respectively.
For this, let $\tau_D^{\varphi_1,\varphi_2 \mp h}$ denote the first
entry time of $(\varPhi^{\varphi_1},\varPhi^{\varphi_2 \mp h})$ into
$D$ for $h>0$ given and fixed. Since $\tau_D^{\varphi_1,\varphi_2
\mp h}$ is optimal for $\hat V(\varphi_1,\varphi_2 \mp h)$ we find
by \eqref{5.8} that
\begin{align} \h{0pc} \label{9.13}
\hat V(\varphi_1,\varphi_2 \mp h) \m \hat V(\varphi_1,\varphi_2) &\ge
\EE_0 \Big[ \int_0^{\tau_D^{\varphi_1,\varphi_2 \mp h}}\!\!\!\!\!\!\!\!
(1 \p \varPhi_t^{\varphi_1}\! \p \varPhi_t^{\varphi_2 \mp h})\, dt +
\hat M \big( \varPhi_{\!\tau_D^{\varphi_1,\varphi_2 \mp h}}^{\varphi_1},
\varPhi_{\!\tau_D^{\varphi_1,\varphi_2 \mp h}}^{\varphi_2 \mp h} \big)
\Big] \\[-3pt] \notag &\h{-8pc}\h{13pt}- \EE_0 \Big[ \int_0^{\tau_D^{\varphi_1,
\varphi_2 \mp h}}\!\!\!\!\!\! (1 \p \varPhi_t^{\varphi_1}\! \p \varPhi_t^
{\varphi_2})\, dt + \hat M \big( \varPhi_{\!\tau_D^{\varphi_1,
\varphi_2 \mp h}}^{\varphi_1},\varPhi_{\!\tau_D^{\varphi_1,\varphi_2
\mp h}}^{\varphi_2} \big) \Big] \\[-1pt] \notag &\h{-8pc}= \EE_0 \Big[ \int_0^
{\tau_D^{\varphi_1,\varphi_2 \mp h}}\!\!\!\!\!\!\!\! (\mp h\:\! \varPhi_
t^2)\, dt \mp c\:\! h\:\! \varPhi_{\!\tau_D^{\varphi_1,\varphi_2 \mp h}}^
2\, I \big( \tau_D^{\varphi_1,\varphi_2 \mp h} \ne \tau_{D_1}^{\varphi_1,
\varphi_2 \mp h} \big) \Big]
\end{align}
for all $h \in (0,h_0)$ with some $h_0>0$ sufficiently small, where
in the final equality we use \eqref{6.4} and \eqref{9.7} similarly
as in \eqref{9.6} above. Dividing both sides of \eqref{9.13} by $\mp
h$ we obtain
\begin{align} \h{0.5pc} \label{9.14}
&\frac{\hat V(\varphi_1,\varphi_2 \m h) \m \hat V(\varphi_1,\varphi_2)}
{-h} \le \EE_0 \Big[ \int_0^{\tau_D^{\varphi_1,\varphi_2-h}}\!\!\!\!\!
\!\!\! \varPhi_t^2\, dt + c\;\! \varPhi_{\!\tau_D^{\varphi_1,\varphi_
2-h}}^2\, I \big( \tau_D^{\varphi_1,\varphi_2-h} \ne \tau_{D_1}^{\varphi_
1,\varphi_2-h} \big) \Big] \\ \label{9.15} &\frac{\hat V(\varphi_1,
\varphi_2 \p h) \m \hat V(\varphi_1,\varphi_2)}{h} \ge \EE_0 \Big[
\int_0^{\tau_D^{\varphi_1,\varphi_2+h}}\!\!\!\!\!\!\!\! \varPhi_t^2\,
dt + c\;\! \varPhi_{\!\tau_D^{\varphi_1,\varphi_2+h}}^2\, I \big(
\tau_D^{\varphi_1,\varphi_2+h} \ne \tau_{D_1}^{\varphi_ 1,\varphi_2+h}
\big) \Big]
\end{align}
for all $h \in (0,h_0)$. Letting $h \downarrow 0$ and using that the
right-hand side in \eqref{9.14} and \eqref{9.15} tends to $c$ by
\eqref{9.1} and the continuity of $\hat V$ we see that \eqref{9.11}
and \eqref{9.12} hold as claimed. Combining \eqref{9.9}+\eqref{9.10}
and \eqref{9.11}+\eqref{9.12} respectively with the fact that $\hat
M_{\varphi_2}(\varphi_1,\varphi_2)=c$ we see that \eqref{9.3} holds
as claimed. This completes the proof. \hfill $\square$

\vspace{16pt}

\textbf{Corollary 15 (\!$C^1$ \!\!regularity\:\!).} \emph{For the
value function $\hat V$ from \eqref{4.18} we have}
\begin{align} \h{4pc} \label{9.16}
&(\varphi_1,\varphi_2) \mapsto \hat V_{\varphi_1}(\varphi_1,\varphi_2)
\;\; \textit{is continuous on}\;\; (0,\infty)^2 \\[1pt] \label{9.17}
&(\varphi_1,\varphi_2) \mapsto \hat V_{\varphi_2}(\varphi_1,\varphi_2)
\;\; \textit{is continuous on}\;\; (0,\infty)^2.
\end{align}

\textbf{Proof.} We have established in Proposition 14 that $\hat V$
is differentiable on $(0,\infty)^2$. By \eqref{8.4} we know that
$\hat V$ is concave on $[0,\infty)^2$. The claims \eqref{9.16} and
\eqref{9.17} then follow from the general fact that concave
differentiable functions are continuously differentiable on open
sets (see e.g.\ \cite[Theorem 2.2.2]{BV}). This completes the proof.
\hfill $\square$

\section{Non-monotonicity of the optimal stopping boundaries}
%%%%%%%%%%%%%%%%%%%%%%%%%%%%%%%%%%%%%%%%%%%%%%%%%%%%%%%%%%%%%

In this section we show that the optimal stopping boundaries $b_0,
b_1, b_2$ are non-monotone as functions of their arguments and prove
the existence of a `belly' which determines their curvature/shape.
In the first part of the proof we introduce the local time of
$\varPhi$ on a fictitious curve which enables us to decompose the
two-dimensional optimal stopping problem into two one-dimensional
optimal stopping problems which can be solved explicitly. In the
second part of the proof we follow the general hint from
\cite[Remark 13]{Pe-3} on establishing the absence of jumps of the
optimal stopping boundaries and make use of Hopf's boundary point
lemma to derive a contradiction with the directional smooth fit. In
view of the symmetry relations \eqref{8.20}+\eqref{8.21} it is
sufficient to focus on the optimal stopping boundary $b_2$ and these
facts then extend to the optimal stopping boundaries $b_0$ and $b_1$
as discussed in Section 8 above.

\vspace{6pt}

1.\ To derive that the optimal stopping set $D_2$ has a `belly' as
displayed on Figure 1 below, we first show that not only the point
$(\beta,0)$ belongs to $D_2$ as derived in Section 7 above but also
a non-trivial vertical segment above $(\beta,0)$ is contained in
$D_2$.

\vspace{12pt}

\textbf{Proposition 16.} \emph{For the stopping set $D_2$ from
\eqref{8.10} we have
\begin{equation} \h{9pc} \label{10.1}
\{ \beta \}\! \times\! [0,\varphi_2] \subseteq D_2
\end{equation}
for some $\varphi_2>0$ small enough.}

\vspace{12pt}

\textbf{Proof.} The idea is to introduce the local time of $\varPhi$
on the line
\begin{equation} \h{4pc} \label{10.2}
c_0' = \{\, (\varphi_1,\varphi_2) \in [0,\infty)^2\; \vert\; \varphi_1
=1\;\; \&\;\; \varphi_2 \in [1,\infty)\,\}
\end{equation}
and decompose the two-dimensional optimal stopping problem
\eqref{4.18} into two one-dimensional optimal stopping problems that
can be solved explicitly.

For this, set $\varphi_1^* := \beta$ throughout and consider the
Lagrange reformulation \eqref{6.5} of the optimal stopping problem
\eqref{4.18} that yields
\begin{equation} \h{-0.5pc} \label{10.3}
\hat V(\varphi_1^*,\varphi_2) \m \hat M(\varphi_1^*,\varphi_2) =
\inf_\tau\;\! \EE_{\varphi_1^*,\varphi_2}^0 \Big[ \int_0^\tau\!
\big( 1 \p \varPhi_t^1 \p \varPhi_t^2 \big)\, dt - \frac{c}{2}\:\!
\Big( \ell_\tau^{c_0}(\varPhi) \p \ell_\tau^{c_1}(\varPhi) \p
\ell_\tau^ {c_2}(\varPhi)\Big)\:\! \Big]
\end{equation}
for $\varphi_2 \in [0,\infty)$ where the infimum is taken over all
stopping times $\tau$ of $\varPhi$. Since the left-hand side of
\eqref{10.3} is non-positive, it is enough to show that the
left-hand side of \eqref{10.3} is non-negative for all $\varphi_2>0$
sufficiently small. For this, adding and subtracting
$\ell_\tau^{c_0'}(\varPhi)$ under the expectation sign in
\eqref{10.3} and noting that
\begin{equation} \h{8pc} \label{10.4}
\ell^{c_0}(\varPhi) + \ell^{c_0'}(\varPhi) = \ell^1(\varPhi^1)
\end{equation}
we find that
\begin{align} \h{2pc} \label{10.5}
\inf_\tau\;\! \EE_{\varphi_1^*,\varphi_2}^0 &\Big[ \int_0^\tau\!
\big( 1 \p \varPhi_t^1 \p \varPhi_t^2 \big)\, dt - \frac{c}{2}\:\!
\Big( \ell_\tau^{c_0}(\varPhi) \p \ell_\tau^{c_1}(\varPhi) \p \ell_
\tau^ {c_2}(\varPhi)\Big)\:\! \Big] \\ \notag &\ge \inf_\tau\;\!
\EE_{\varphi_1^*,\varphi_2}^0 \Big[ \int_0^\tau\! \big( 1 \p \varPhi_
t^1 \big)\, dt - \frac{c}{2}\;\! \ell_\tau^1(\varPhi^1) \:\! \Big]
\\ \notag &\h{13pt}+ \inf_\tau\;\! \EE_{\varphi_1^*,\varphi_2}^0
\Big[ \int_0^\tau\! \varPhi_t^2\, dt + \frac{c}{2}\;\! \ell_\tau^
{c_0'}(\varPhi) - \frac{c}{2}\:\! \Big( \ell_\tau^{c_1}(\varPhi)
\p \ell_\tau^ {c_2}(\varPhi) \Big)\:\! \Big] \\ \notag &= \inf_
\tau\;\! \EE_{\varphi_1^*,\varphi_2}^0 \Big[ \int_0^\tau\! \big(
1 \p \varPhi_t^1 \big)\, dt + c\:\! \big( 1\! \wedge\! \varPhi_
\tau^1 \big)\:\! \Big] - c \:\! (1\! \wedge\! \varphi_1^*) \\
\notag &\h{13pt}+ \inf_\tau\;\! \EE_{\varphi_1^*,\varphi_2}^0
\Big[ \int_0^\tau\! \varPhi_t^2\, dt + \tilde M(\varPhi_\tau^1,
\varPhi_\tau^2)\:\! \Big] - \tilde M(\varphi_1^*,\varphi_2) \\
\notag &\ge \inf_\tau\;\! \EE_{\varphi_1^*,\varphi_2}^0 \Big[
\int_0^\tau\! \varPhi_t^2\, dt + c\:\! \big( 1\! \wedge\!
\varPhi_\tau^2 \big)\:\! \Big] - c \:\! (1\! \wedge\! \varphi_2)
\end{align}
for $\varphi_2 \in [0,1]$ where in the equality we use the
It\^o-Tanaka formula (cf.\ \cite[p.\ 223]{RY}) applied to $c\:\! (1
\wedge \varPhi^1)$, and the change-of-variable formula with local
time on surfaces \cite[Theorem 2.1]{Pe-2} applied to $\tilde
M(\varPhi^1,\varPhi^2)$ similarly to \eqref{6.9} above with
\begin{equation} \h{8pc} \label{10.6}
\tilde M := c\:\! \big[ (1 \vee \varphi_1) \wedge \varphi_2 \big]
\end{equation}
for $(\varphi_1,\varphi_2) \in [0,\infty)^2$\!, both combined with
the optional sampling theorem upon using that $\varPhi^1$ and
$\varPhi^2$ are martingales under $\PP_{\!0}$. In the final
inequality of \eqref{10.5} we use that $\varphi_1^* = \beta$ is an
optimal stopping point in the one-dimensional optimal stopping
problem for $\varPhi^1$ as established in Section 7 above as well as
that $\tilde M(\varphi_1^*,\varphi_2) = c\:\! [ (1 \vee \varphi_1^*)
\wedge \varphi_2 ] =  c\:\! (\varphi_1^* \wedge \varphi_2) = c\:\!
\varphi_2 = c\:\! (1 \wedge \varphi_2)$ for $\varphi_2 \in [0,1]$ as
claimed.

%%%%%%%%%%%%%%%%%%%%%%%%%%%%%%%%%%%%%%%%%%%%%%%%%%%%%%%%%%%%%%%%%%%%%%%%%%%%%%%
%%% Figure 1 %%%
%%%%%%%%%%%%%%%%%%%%%%%%%%%%%%%%%%%%%%%%%%%%%%%%%%%%%%%%%%%%%%%%%%%%%%%%%%%%%%%

\begin{figure}[!t]
\begin{center}
\includegraphics[scale=0.7]{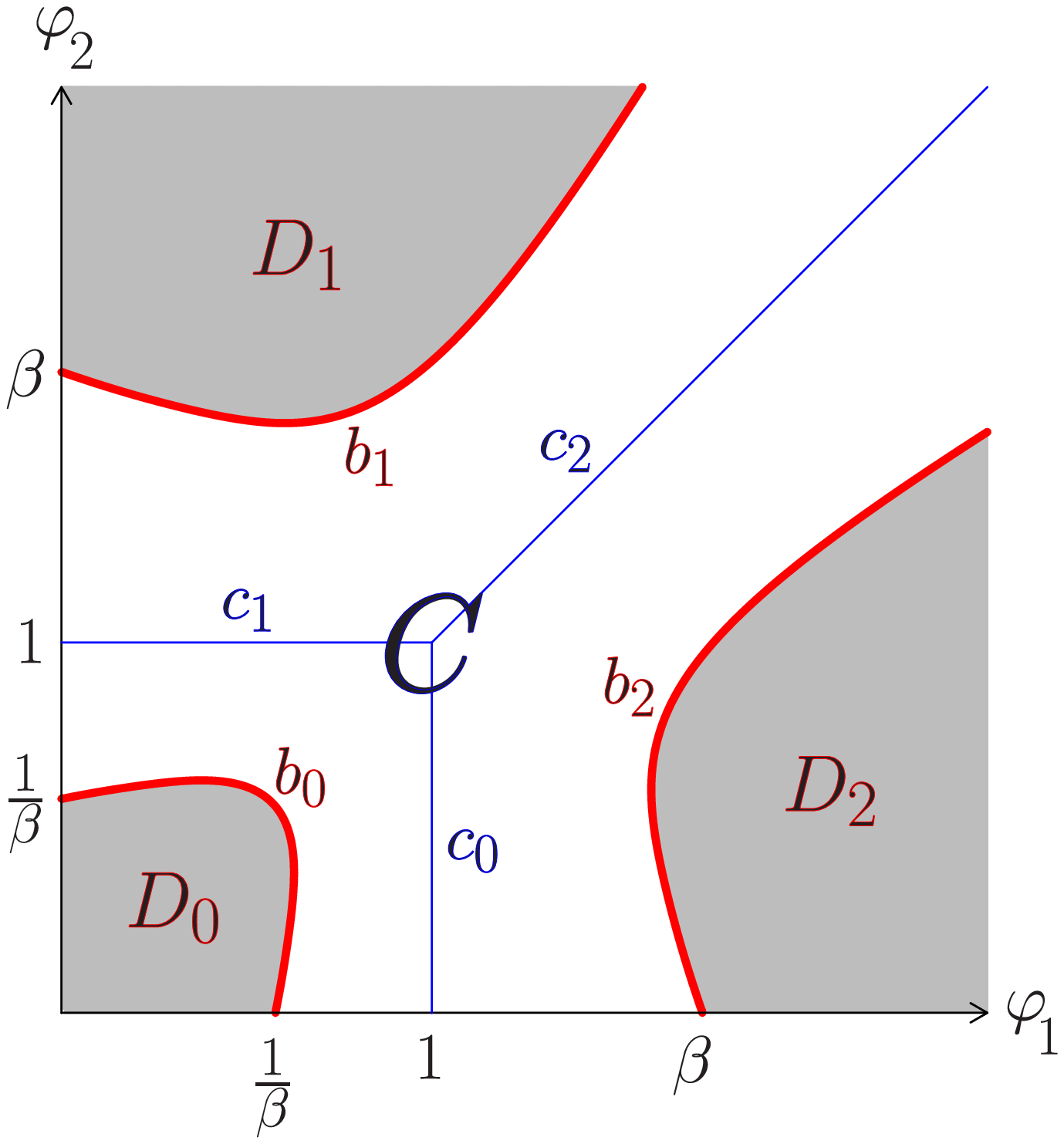}
\end{center}

{\par \leftskip=1.8cm \rightskip=1.6cm \small \noindent \vspace{-16pt}

\textbf{Figure 1.} Location of the continuation set $C$, the
stopping sets $D_0, D_1, D_2$, and the optimal stopping boundaries
$b_0, b_1, b_2$, recalling that $b_1(\varphi_1) = b_2(\varphi_2)$
for $\varphi_1 = \varphi_2$ in $[0,\infty)$.

\par} \vspace{0pt}

\end{figure}

%%%%%%%%%%%%%%%%%%%%%%%%%%%%%%%%%%%%%%%%%%%%%%%%%%%%%%%%%%%%%%%%%%%%%%%%%%%%%%%
%%% End Figure 2 %%%
%%%%%%%%%%%%%%%%%%%%%%%%%%%%%%%%%%%%%%%%%%%%%%%%%%%%%%%%%%%%%%%%%%%%%%%%%%%%%%%

Motivated by the right-hand side in \eqref{10.5} above, consider the
optimal stopping problem
\begin{equation} \h{6pc} \label{10.7}
\tilde V (\varphi) = \inf_\tau\;\! \EE_{\varphi}^0 \Big[ \int_0^\tau
\varPhi_t\, dt + c\:\! \big( 1\! \wedge\! \varPhi_\tau
\big) \:\! \Big]
\end{equation}
for $\varphi \in [0, \infty)$ with $\PP_{\!\varphi}^0(\varPhi_0\!
=\! \varphi) = 1$ where the process $\varPhi$ and its infinitesimal
generator $\LL_\varPhi$ are given by \eqref{7.2}+\eqref{7.3} and
\eqref{7.4} above, and the infimum in \eqref{10.7} is taken over all
stopping times $\tau$ of $\varPhi$. The optimal stopping problem
\eqref{10.7} is similar to the optimal stopping problem \eqref{7.1}
and we can use similar arguments to tackle it. Denoting the loss
function in \eqref{10.7} by $\tilde M(\varphi) = c\:\! (1\! \wedge\!
\varphi)$ for $\varphi \in [0,\infty)$ it follows that the
free-boundary problem now reads
\begin{align} \h{4pc} \label{10.8}
&\LL_\varPhi \tilde V(\varphi) = -\varphi\;\; \text{for}\;\; \varphi
\in (\tilde \varphi_0^*,\tilde \varphi_1^*) \\[1pt] \label{10.9} &\tilde
V(\tilde \varphi_i^*) = \tilde M(\tilde \varphi_i^*)\;\; \text{for}
\;\; i=0,1\;\; \text{(instantaneous stopping)} \\[1pt] \label{10.10}
&\tilde V'(\tilde \varphi_i^*) = \tilde M'(\tilde \varphi_i^*)
\;\; \text{for}\;\; i=0,1\;\; \text{(smooth fit)}
\end{align}
where $0 < \tilde \varphi_0^* < 1 < \tilde \varphi_1^* < \infty$ are
the optimal stopping/boundary points to be found and we have $\tilde
V(\varphi) = \tilde M(\varphi)$ for $\varphi \in [0,\tilde
\varphi_0^*) \cup (\tilde \varphi_1^*,\infty)$ as well (in addition
to \eqref{10.9} above).

The general solution to the ordinary differential equation
\eqref{10.8} is given by
\begin{equation} \h{6pc} \label{10.11}
\hat V(\varphi) = \tilde A\:\! \varphi + \tilde B + \frac{1}{\mu^2}\;\!
\varphi\:\! (1 \m \log \varphi)
\end{equation}
for $\varphi > 0$ where $\tilde A$ and $\tilde B$ are two
undetermined real constants. Boundary conditions \eqref{10.9} and
\eqref{10.10} then read as follows
\begin{align} \h{6pc} \label{10.12}
&\tilde A\:\! \tilde \varphi_0^* + \tilde B + \frac{1}{\mu^2}\;\!
\tilde \varphi_0^*\:\! (1 \m \log \tilde \varphi_0^*) = c\:\! \varphi_0^*
\\ \label{10.13} &\tilde A\:\! \tilde \varphi_1^* + \tilde B +
\frac{1}{\mu^2}\;\! \tilde \varphi_1^*\:\! (1 \m \log \tilde \varphi_1^*)
= c \\ \label{10.14} &\tilde A - \frac{1}{\mu^2}\:\! \log \tilde
\varphi_0^* = c \\ \label{10.15} &\tilde A - \frac{1}{\mu^2}\:\!
\log \tilde \varphi_1^* = 0\, .
\end{align}
It is a matter of routine to verify that the unique solution to
\eqref{10.12}-\eqref{10.15} is given by
\begin{align} \h{6pc} \label{10.16}
&\tilde \varphi_0^* = \frac{c \mu^2}{e^{c \mu^2} \m 1} \quad \& \quad
\tilde \varphi_1^* = \frac{c \mu^2 e^{c \mu^2}}{e^{c \mu^2} \m 1} \\
\label{10.17} &\tilde A^* = c \p  \frac{1}{\mu^2}\:\! \log \tilde
\varphi_0^* \quad \& \quad \tilde B^* = -\frac{1}{\mu^2}\;\! \tilde
\varphi_0^*\, .
\end{align}
Note that $\tilde \varphi_0^* \in (0,1)$ and $\tilde \varphi_1^* \in
(1,\infty)$ as needed. Inserting $\tilde A^*$ and $\tilde B^*$ from
\eqref{10.17} to \eqref{10.11} we obtain a candidate value function
$\tilde V^*$ for the optimal stopping problem \eqref{10.7}. Applying
the It\^o-Tanaka formula (cf.\ \cite[p.\ 223]{RY}) to $\tilde V^*$
composed with $\varPhi$, which reduces to It\^o's formula due to
smooth fit \eqref{10.10}, and making use of the optional sampling
theorem, it is easily verified that $\tilde V^*$ coincides with the
value function $\tilde V$ from \eqref{10.7} and the optimal stopping
time (at which the infimum in \eqref{10.7} is attained) is given by
\begin{equation} \h{6pc} \label{10.18}
\tau_* = \inf\, \{\, t \ge 0\; \vert\; \varPhi_t \notin (\tilde
\varphi_0^*, \tilde \varphi_1^*)\, \}
\end{equation}
where $\tilde \varphi_0^*$ and $\tilde \varphi_1^*$ are given by
\eqref{10.16} above. This in particular shows that the interval
$[0,\tilde \varphi_0^*]$ is contained in the stopping set of the
optimal stopping problem \eqref{10.7}. Translating this conclusion
to the right-hand side of \eqref{10.5} above we see that its value
equals zero whenever $\varphi_2$ belongs to $[0,\tilde
\varphi_0^*]$. It follows therefore from \eqref{10.3} and
\eqref{10.5} that $\varphi_2$ in \eqref{10.1} can be taken to be
equal to $\tilde \varphi_0^* = c \mu^2\!/(e^{c \mu^2} \m 1)$ and the
proof is complete. \hfill $\square$

\vspace{12pt}

2.\ We now show that the `belly' of the optimal stopping set $D_2$
is not flat but curved (see Figure 1 above). For this, suppose that
this is not the case. Then $[a,b)\! \times\! [c,d] \subseteq C$ with
$\{ b \}\! \times\! [c,d] \subseteq D_2$ and $\{ b \}\! \times\!
(d,b] \subseteq C$ for some $a<b$ with $[a,b] \subseteq [1,\beta]$
and some $c<d$ with $[c,d] \subseteq [0,a]$. The initial claim is
then a direct consequence of the following fact.

\vspace{12pt}

\textbf{Proposition 17.} \emph{If the `belly' of the optimal
stopping set $D_2$ would be flat as described above, then the
horizontal smooth fit condition \eqref{9.2} would fail on $\{ b \}\!
\times\! [c,d] \subseteq \partial C \cap D_2$.}

\vspace{12pt}

\textbf{Proof.} Suppose that the `belly' of the optimal stopping set
$D_2$ is flat as described above. Set $R^0 = (a,b)\! \times\! (c,d)$
and $R^1 = (a,b]\! \times\! (c,d)$ with $R = [a,b]\! \times\!
[c,d]$. Recalling the Lagrange reformulation \eqref{6.5} of the
optimal stopping problem \eqref{4.18}, and arguing as in the proof
of Theorem 12 in \cite{Pe-3}, we find that the value function $\hat
V$ from \eqref{4.18} solves the equation
\begin{equation} \h{9.5pc} \label{10.19}
\LL_\varPhi \hat V = - \hat H
\end{equation}
on $R^0$ and belongs to $C^4(R^1)$ where $\LL_\varPhi$ is given by
\eqref{5.9} above and we set
\begin{equation} \h{8pc} \label{10.20}
\hat H(\varphi_1,\varphi_2) = 1 \p \varphi_1 \p \varphi_2
\end{equation}
for $(\varphi_1,\varphi_2) \in [0,\infty)^2$. Differentiating both
sides of \eqref{10.19} with respect to $\varphi_2$ and defining the
differential operator $\tilde \LL$ by setting
\begin{equation} \h{2pc} \label{10.21}
\tilde \LL = \varphi_1^2\;\! \partial_{\varphi_1 \varphi_1}^2 +
\varphi_1 \varphi_2\;\! \partial_{\varphi_1 \varphi_2}^2
+ \varphi_2^2\;\! \partial_{\varphi_2 \varphi_2}^2 +
2 \varphi_1\;\! \partial_{\varphi_1} + 4 \varphi_2\;\!
\partial_{\varphi_2} + 2
\end{equation}
we find that $\hat V_{\varphi_2 \varphi_2}$ solves the equation
\begin{equation} \h{9.5pc} \label{10.22}
\tilde \LL\:\! \hat V_{\varphi_2 \varphi_2} = 0
\end{equation}
on $R^0$. We will now complete the proof in two steps as follows.

\vspace{6pt}

1.\ We claim that the strict inequality holds
\begin{equation} \h{9pc} \label{10.23}
\hat V_{\varphi_2 \varphi_2}(\varphi_1,\varphi_2) < 0
\end{equation}
for all $(\varphi_1,\varphi_2) \in R^0$. For this, suppose that
\eqref{10.23} fails for some $(\varphi_1,\varphi_2) \in R^0$.
Recalling that $\{ b \}\! \times\! (d,b] \subseteq C$, consider the
ball $b(z,r)$ with centre at $z:=(b,d)$ and radius $r>0$ small
enough so that $b(z,r) \subseteq \Delta_2$, where $\Delta_2$ is
defined following \eqref{6.4} above. Enlarge $R^0$ by setting
$\tilde R^0 := R^0 \cup (b(z,r) \cap C)$ and note that the same
arguments as above show that the equations \eqref{10.19} and
\eqref{10.21} hold on $\tilde R^0$ too. Since the coefficients of
$\tilde \LL$ are continuous and the set $\tilde R^0$ is bounded we
can conclude that $\tilde \LL$ is uniformly elliptic on $\tilde R^0$
(cf.\ \cite[p.\ 31]{GT}). The hypothesis that \eqref{10.23} fails
for some $(\varphi_1,\varphi_2) \in R^0$\!, combined with the fact
that $\hat V_{\varphi_2 \varphi_2} \le 0$ on $\tilde R^0$ by
\eqref{8.4} above, implies that $\hat V_{\varphi_2
\varphi_2}(\varphi_1,\varphi_2 )=0$ so that $\hat V_{\varphi_2
\varphi_2}$ attains its maximum in the interior of $\tilde R^0$
(i.e.\ not at its boundary alone). Hence by the strong maximum
principle for elliptic equations (see Theorem 3.5 in \cite[p.\
35]{GT} and the second sentence following its proof) we can conclude
that $\hat V_{\varphi_2 \varphi_2} = 0$ on the entire $\tilde R^0$.
This in particular means that $\varphi_2 \mapsto \hat
V(b,\varphi_2)$ is linear on $[d,d \p r]$. Since $\varphi_2 \mapsto
\hat V(b,\varphi_2) = c\:\! (1 \p \varphi_2)$ is linear on $[c,d]$
as well, and the vertical smooth fit \eqref{9.3} holds at $z=(b,d)$,
it follows that $\hat V(b,\varphi_2) = c\:\! (1 \p \varphi_2)$ for
all $\varphi_2 \in [c,d \p r]$ so that $\{ b \}\! \times\! (d,d \p
r] \subseteq D$ which is a contradiction. This establishes that
\eqref{10.23} is satisfied as claimed.

\vspace{6pt}

2.\ Fix any point $e$ in $(c,d)$ and note that $\hat V_{\varphi_2
\varphi_2}(b,e) = 0$ since $\hat V \in C^4(R^1) \subseteq C^2(R^1)$
and $\hat V(b,\varphi_2) = c\: (1 \p \varphi_2)$ for $\varphi_2 \in
[c,d]$. Hence we see that \eqref{10.23} reads as $\hat V_{\varphi_2
\varphi_2}(\varphi_1,\varphi_2) < \hat V_{\varphi_2 \varphi_2}(b,e)$
for all $(\varphi_1,\varphi_2) \in R^0$. Moreover, we know that
$\tilde \LL$ is uniformly elliptic and $\tilde \LL\:\! \hat
V_{\varphi_2 \varphi_2} \ge 0$ holds on $R^0$ by \eqref{10.22}
above. Finally, it is evident that $R^0$ satisfies an interior
sphere condition at $z := (b,e) \in \partial R^0$ (\;\!i.e.\ there
exist $w \in R^0$ and $r>0$ such that $b(w,r) \subseteq R^0$ and $z
\in
\partial (b(w,r))$). These facts show that Hopf's boundary point
lemma for elliptic equations (see \cite[Lemma 3.4 p.\ 34]{GT}) is
applicable and thus the outer normal derivative of $\hat
V_{\varphi_2 \varphi_2}$ at $z=(b,e)$ must be strictly positive. In
other words, we have
\begin{equation} \h{8pc} \label{10.24}
\big( \hat V_{\varphi_2 \varphi_2} \big)_{\!\varphi_1}\!(b,e) > 0\, .
\end{equation}
This conclusion shows that the horizontal smooth fit condition
\eqref{9.2} cannot hold on $\{ b \}\! \times\! [c,d]$ as claimed,
since otherwise we would have $(\hat V_{\varphi_2
\varphi_2})_{\varphi_1}(b,e) = (\hat V_{\varphi_1})_{\varphi_2
\varphi_2}(b,e) = 0$ due to $\hat V \in C^4(R^1) \subseteq
C^3(R^1)$, and the proof is complete. \hfill $\square$

\section{Free-boundary problem}
%%%%%%%%%%%%%%%%%%%%%%%%%%%%%%%

In this section we derive a free-boundary problem that stands in
one-to-one correspondence with the optimal stopping problem
\eqref{4.18} and establish the fact that the value function $\hat V$
and the optimal stopping boundary $\partial C$ solve the
free-boundary problem uniquely. These considerations will be
continued in the next section.

\vspace{6pt}

1.\ Consider the optimal stopping problem \eqref{4.18} where the
strong Markov/Feller process $\varPhi = (\varPhi^1,\varPhi^2)$
solves the system of stochastic differential equations
\eqref{5.6}+\eqref{5.7}. Recalling that the infinitesimal generator
$\LL_\varPhi$ of $\varPhi$ is given by \eqref{5.9} above, and
relying on other properties of $\hat V$ and $\partial C$ derived in
Section 8 above, we are naturally led to formulate the following
free-boundary problem for finding $\hat V$ and $\partial C$:
\begin{align} \h{5pc} \label{11.1}
&\LL_\varPhi \hat V = - \hat H\;\; \text{on}\;\; C \\[1pt]
\label{11.2} &\hat V = \hat M\;\; \text{on}\;\; D\;\; \text{
(instantaneous stopping)} \\[1pt]  \label{11.3} &\hat V_{
\varphi_i} = \hat M_{\varphi_i} \;\; \text{on}\;\; \partial
C\;\; \text{for}\;\; i=1,2\;\; \text{(smooth fit)}
\end{align}
where $\hat H$ is given by \eqref{10.20} above and $\hat M$ is given
by \eqref{4.13} above. The continuation set $C$ and the stopping set
$D$ are formally defined by \eqref{8.1} and \eqref{8.2}
respectively. We know from the results of Section 8 that the optimal
stopping boundary $\partial C$ can be fully described by means of
the functions $b_0$ and $b_1$ defined in Section 8 above via the
equivalence $(\varphi_1,\varphi_2) \in \partial C$ if and only if
either $(\varphi_1,\varphi_2) \in \partial C \cap D_0$ and
$\varphi_i = b_0(\varphi_j)$ when $\varphi_i \ge \varphi_j$ for $i
\ne j \in \{1,2\}$ or $(\varphi_1,\varphi_2) \in \partial C \cap
D_i$ and $\varphi_j = b_1(\varphi_i)$ for $i \ne j \in \{1,2\}$
where $D_0, D_1, D_2$ are given by \eqref{8.8}-\eqref{8.10} above
(see Figure 1 above). Clearly the global condition \eqref{11.2} can
be replaced by the local condition $\hat V = \hat M$ on $\partial C$
so that the free-boundary problem \eqref{11.1}-\eqref{11.3} needs to
be considered on the closure of $C$ only (\:\!extending $\hat V$ to
the rest of $D$ as $\hat M$\!).

\vspace{6pt}

2.\ To formulate the existence and uniqueness result for the
free-boundary problem \eqref{11.1}-\eqref{11.3} we let $\cal C$
denote the class of functions $(F;a_0,a_1)$ such that
\begin{align} \h{0pc} \label{11.4}
&F\;\; \text{is concave and continuous on}\;\; [0,\infty)^2\;\;
\text{and belongs to}\;\; C^1((0,\infty)^2) \cap C^2(C_{a_0,a_1})
\\[4pt] \label{11.5} &a_0\;\; \text{is concave and continuous
on}\;\; [0,\delta]\;\; \text{with}\;\; a_0(0) = 1/\beta\;\;
\&\;\; a_0(\delta)=\delta\;\; \\[-1pt] \notag &\text{and}\;\;
\varphi < a_0(\varphi) < 1\;\; \text{for}\;\; \varphi \in (0,\delta)
\;\; \text{with some} \;\; \delta \in (0,1) \\[4pt] \label{11.6}
&a_1\;\; \text{is convex and continuous on}\;\; [0,\infty)\;\;
\text{with}\;\; a_1(0) = \beta\;\; \&\;\; a_1(\infty) = \infty
\;\; \\[-1pt] \notag &\text{and}\;\; a_1(\varphi) > 1 \vee
\varphi\;\; \text{for}\;\; \varphi \in [0,\infty)
\end{align}
where $C_{a_0,a_1} := \{\, (\varphi_1,\varphi_2) \in [0,\infty)^2\;
\vert\; a_0(\varphi_i) < \varphi_j < a_1(\varphi_i)\;\;
\text{when}\;\; \varphi_i \le \varphi_j\;\; \text{and}\;\; \varphi_i
\in [0,\delta]\;\; \text{or}\;\; \varphi_i \le \varphi_j <
a_1(\varphi_i)\;\; \text{when}\;\; \varphi_i \le \varphi_j\;\;
\text{and}\;\; \varphi_i \in (\delta,\infty)\;\; \text{for}\;\; i
\ne j \in \{1,2\}\;\; \text{and some}\;\; \delta \in (0,1)\, \}$ is
the open set surrounded by $a_0$ and $a_1$ (applied twice).

\vspace{12pt}

\textbf{Theorem 18.} \emph{The free-boundary problem
\eqref{11.1}-\eqref{11.3} has a unique solution $(\hat V;b_0,b_1)$
in the class $\cal C$ where $\hat V$ is given by \eqref{4.18} while
$b_0$ and $b_1$ are defined in Section 8 above.}

\vspace{12pt}

\textbf{Proof.} Combining the results of Proposition 5 and Corollary
15 with the arguments leading to \eqref{10.19} above, we see that
the value function $\hat V$ from \eqref{4.18} satisfies \eqref{11.4}
and solves the boundary value problem \eqref{11.1}-\eqref{11.3} with
$\partial C$ described by $b_0$ and $b_1$ from Section 8 as recalled
above. Moreover, combining the results of Propositions 8-12 we see
that $b_0$ and $b_1$ satisfy \eqref{11.5} and \eqref{11.6}
respectively. This shows that $(\hat V;b_0,b_1)$ solves the
free-boundary problem \eqref{11.1}-\eqref{11.3} in the class $\cal
C$ as claimed. To derive uniqueness of the solution we will first
see in the next section that any solution $(F;a_0,a_1)$ to the
free-boundary problem \eqref{11.1}-\eqref{11.3} in the class $\cal
C$ admits a closed triple-integral representation of $F$ expressed
in terms of $a_0$ and $a_1$, which in turn solve a coupled system of
nonlinear Fredholm integral equations, and we will see that this
system cannot have other solutions satisfying the specified
properties. Drawing these facts together we can conclude that there
cannot exist more than one solution to the free-boundary problem
\eqref{11.1}-\eqref{11.3} in the class $\cal C$ as claimed. \hfill
$\square$

\section{Nonlinear integral equations}
%%%%%%%%%%%%%%%%%%%%%%%%%%%%%%%%%%%%%%

In this section we show that the optimal stopping boundaries $b_0$
and $b_1$ can be characterised as the unique solution to a coupled
system of nonlinear Fredholm integral equations (recall that $b_2$
coincides with $b_1$ in terms of its functional rule). This also
yields a closed triple-integral representation of the value function
$\hat V$ expressed in terms of the optimal stopping boundaries $b_0$
and $b_1$. As a consequence of the existence and uniqueness result
for the coupled system of nonlinear Fredholm integral equations we
also obtain uniqueness of the solution to the free-boundary problem
\eqref{11.1}-\eqref{11.3} as explained in the proof of Theorem 18
above. Finally, collecting the results derived throughout the paper
we conclude our exposition at the end of this section by disclosing
the solution to the initial problem.

\vspace{6pt}

1.\ To formulate the theorem below, let $p$ denote the transition
probability density function of the (time-homogeneous) Markov
process $\varPhi = (\varPhi^1,\varPhi^2)$ under $\PP_{\!0}$ in the
sense that
\begin{equation} \h{3pc} \label{12.1}
\PP_{\!\varphi_1,\varphi_2}^0 \big( \varPhi_t^1 \le \psi_1,
\varPhi_t^2 \le \psi_2 \big) = \int_0^{\psi_1}\!\!\! \int_0^{\psi_2}
p(t;\varphi_1,\varphi_2;\eta_1,\eta_2)\, d\eta_1 d\eta_2
\end{equation}
for $t>0$ with $(\varphi_1,\varphi_2)$ and $(\psi_1,\psi_2)$ in
$[0,\infty)^2$. A lengthy but straightforward calculation based on
\eqref{5.8} shows that
\begin{align} \h{-0.4pc} \label{12.2}
p(t;\varphi_1,\varphi_2;\psi_1,\psi_2) = \frac{1}{2 \pi\:\! \sqrt{3}
\, \mu^2\, t\, \psi_1\:\! \psi_2} \exp\! \bigg[\! &- \frac{1}{3}
\bigg( \mu^2 t + \log\! \Big( \frac{\psi_1\:\! \psi_2}{\varphi_1\:\!
\varphi_2} \Big) \\ \notag &+ \frac{1}{\mu^2 t} \Big[ \log^2\! \Big(
\frac{\psi_1\:\! \varphi_2}{\psi_2\:\! \varphi_1} \Big) + \log\!
\Big( \frac{\psi_1}{\varphi_1} \Big) \log\! \Big( \frac{\psi_2}
{\varphi_2} \Big) \Big] \bigg) \bigg]
\end{align}
for $t>0$ with $(\varphi_1,\varphi_2)$ and $(\psi_1,\psi_2)$ in
$[0,\infty)^2$. Recalling from Section 8 above that the functions
$b_0$ and $b_1$ are sufficient to describe the entire boundary of
the continuation set, we can then evaluate the expression of
interest in the theorem below as follows
\begin{align} \h{3pc} \label{12.3}
&\EE_{\varphi_1,\varphi_2}^0 \Big[ \int_0^\infty \hat H(\varPhi_s^1,
\varPhi_s^2)\, I \big( (\varPhi_s^1,\varPhi_s^2)\! \in\! C \big)\,
ds\;\! \Big] \\ \notag &= \int_0^\infty\!\! \bigg( \int_0^\infty\!\!
\int_{b_0(\psi_1) \vee \psi_1}^{b_1(\psi_1)}\!\!\! \hat H(\psi_1,
\psi_2)\, p(s;\varphi_1,\varphi_2;\psi_1,\psi_2)\: d \psi_2\;\!
d \psi_1 \\ \notag &\h{42pt}+ \int_0^\infty\!\! \int_{b_0(\psi_2)
\vee \psi_2}^{b_1(\psi_2)}\!\!\! \hat H(\psi_1,\psi_2)\, p(s;
\varphi_1,\varphi_2;\psi_1,\psi_2)\: d \psi_1\;\! d \psi_2 \bigg)
\;\! ds \\ \notag &= 2 \int_0^\infty\!\! \int_0^\infty\!\!
\int_{b_0(\psi_1) \vee \psi_1}^{b_1(\psi_1)}\!\!\! \hat H(\psi_1,
\psi_2)\, p(s;\varphi_1,\varphi_2;\psi_1,\psi_2)\: d \psi_2\, d
\psi_1\;\! ds
\end{align}
for $(\varphi_1,\varphi_2) \in [0,\infty)^2$ where the final
equality follows by symmetry relative to the main diagonal in
$[0,\infty)^2$ and we recall that $\hat H$ is defined in
\eqref{10.20} above.

\vspace{12pt}

\textbf{Theorem 19 (Existence and uniqueness).} \emph{The optimal stopping
boundaries $b_0$ and $b_1$ in the problem \eqref{4.18} can be
characterised as the unique solution to the coupled system of
nonlinear Fredholm integral equations
\begin{align} \h{2pc} \label{12.4}
&\varphi_1 + b_0(\varphi_1) = \frac{2}{c} \int_0^\infty\!\! \int_0
^\infty\!\! \int_{b_0(\psi_1) \vee \psi_1}^{b_1(\psi_1)}\!\!\! \hat
H(\psi_1,\psi_2)\, p(s;\varphi_1,\varphi_2;\psi_1,\psi_2)\: d \psi_2
\, d \psi_1\;\! ds \\ \label{12.5} &1 + \varphi_1 = \frac{2}{c} \int
_0^\infty\!\! \int_0^\infty\!\! \int_{b_0(\psi_1) \vee \psi_1}^{b_1
(\psi_1)}\!\!\! \hat H(\psi_1,\psi_2)\, p(s;\varphi_1,\varphi_2;
\psi_1,\psi_2)\: d \psi_2\, d \psi_1\;\! ds
\end{align}
in the class of functions $a_0$ and $a_1$ satisfying \eqref{11.5}
and \eqref{11.6} respectively, where $\varphi_1$ in \eqref{12.4}
belongs to $[0,\gamma]$ with $b_0(\gamma) = \gamma$ for some $\gamma
\in (0,1)$ and $\varphi_1$ in \eqref{12.5} belongs to $[0,\infty)$.
The value function $\hat V$ in the problem \eqref{4.18} admits the
following representation
\begin{equation} \h{2pc} \label{12.6}
\hat V(\varphi_1,\varphi_2) = 2 \int_0^\infty\!\! \int_0^\infty\!\!
\int_{b_0(\psi_1) \vee \psi_1}^{b_1(\psi_1)}\!\!\! \hat H(\psi_1,
\psi_2)\, p(s;\varphi_1,\varphi_2; \psi_1,\psi_2)\: d \psi_2\, d
\psi_1\;\! ds
\end{equation}
for $(\varphi_1,\varphi_2) \in [0,\infty)^2$. The optimal stopping
time in the problem \eqref{4.18} is given by
\begin{equation} \h{-0.5pc} \label{12.7}
\tau_{b_0,b_1} = \inf\, \{\, t \ge 0\; \vert\; \varPhi_t^i \ge \varPhi_t^j
\;\; \text{with}\;\; \varPhi_t^i \le b_0(\varPhi_t^j)\;\; \text{or}\;\;
\varPhi_t^i \ge b_1(\varPhi_t^j)\;\; \text{for}\;\; i \ne j\;\; \text{in}
\;\; \{1,2\}\, \}
\end{equation}
under $\PP_{\!\varphi_1,\varphi_2}^0$ with $(\varphi_1,\varphi_2)
\in [0,\infty)^2$ given and fixed (see Figure 1 above).}

\vspace{12pt}

\textbf{Proof.} (I) \emph{Existence}. We first show that the value
function $\hat V$ in the problem \eqref{4.18} admits the
representation \eqref{12.6} and that the optimal stopping boundaries
$b_0$ and $b_1$ solve the system \eqref{12.4}+\eqref{12.5}.
Recalling that $b_0$ and $b_1$ satisfy the properties \eqref{11.5}
and \eqref{11.6} this will establish the existence of a solution to
the system \eqref{12.4}+\eqref{12.5}.

For this, recall that by \eqref{8.4} in Proposition 5 we know that
$\hat V$ is concave and from Corollary 15 we know that $\hat V$ is
globally $C^1$ on $(0,\infty)^2$. These properties however are
generally insufficient to apply a known extension of It\^o's formula
to $\hat V(\varPhi^1,\varPhi^2)$ due to not knowing the size of the
second partial derivatives $\hat V_{\varphi_1,\varphi_1}$, $\hat
V_{\varphi_1,\varphi_2}$, $\hat V_{\varphi_2,\varphi_2}$ close to
the optimal stopping boundaries. Note that we know that the optimal
stopping boundaries are convex/concave, however, this is generally
insufficient to derive a local boundedness of the second partial
derivatives close to the optimal stopping boundaries (without having
their smoothness) using the generally theory of elliptic PDEs (see
\cite{GT}). A semimartingale decomposition of $\hat
V(\varPhi^1,\varPhi^2)$ obtained by It\^o's formula is useful
because it leads to Dynkin's formula (upon localising, taking
expectations, and passing to the limit) which in turn yields the
representation \eqref{12.6}. We will show in the proof below that
Dynkin's formula can be derived without appealing to It\^o's formula
and/or without formally verifying that the second partial
derivatives are locally bounded close to the optimal stopping
boundaries. This will be accomplished in several steps below by
exploiting the underlying convexity/concavity in the problem
\eqref{4.18} combined with the fact that the expectation of the
running local time of $(\varPhi^1,\varPhi^2)$ on the (approximating)
optimal stopping boundaries remains uniformly bounded as the time
tends to infinity (recall that $(\varPhi^1,\varPhi^2)$ itself
converges to zero so that this is rather intuitive).

\vspace{6pt}

1.\ We begin by localising the process $\varPhi =
(\varPhi^1,\varPhi^2)$. For this, let $N \ge 1$ be given and fixed
(large) and consider the first exit time of $\varPhi$ from the
square $[1/N,N]^2$ given by
\begin{equation} \h{7pc} \label{12.8}
\rho_N = \inf\, \{\, t \ge 0\; \vert\; \varPhi_t \notin [1/N,N]^2
\, \}\, .
\end{equation}
Let $\varPhi^{\rho_N} = (\varPhi_{t \wedge \rho_N})_{t \ge 0}$
denote the process $\varPhi$ stopped at $\rho_N$. Clearly the
process $\varPhi^{\rho_N}$ stays in the square $[1/N,N]^2$ all the
time while both $\hat V$ and $\hat V_{\varphi_i}$ are continuous and
thus bounded on $[1/N,N]^2$ for $i=1,2$. As we have not established
that $\hat V_{\varphi_1,\varphi_1}$, $\hat V_{\varphi_1,\varphi_2}$,
$\hat V_{\varphi_2,\varphi_2}$ are locally bounded close to the
optimal stopping boundaries, we proceed by modifying the value
function $\hat V$ within the continuation set $C$ close to its
boundary.

\vspace{6pt}

2.\ For $n \ge 1$ given and fixed (large) define the sets $C_n :=
\{\, (\varphi_1,\varphi_2) \in [0,\infty)^2\; \vert\; \hat
V(\varphi_1,\varphi_2) < \hat M(\varphi_1,\varphi_2) \m 1/n\, \}$
and $D_n := \{\, (\varphi_1,\varphi_2) \in [0,\infty)^2\; \vert\;
\hat V(\varphi_1,\varphi_2) \ge \hat M(\varphi_1,\varphi_2) \m 1/n\,
\}$. Set $C_i^n := C_n \cap \Delta_i$ and $D_i^n := D_n \cap
\Delta_i$ where $\Delta_i$ are defined following \eqref{6.4} above
for $i=0,1,2$. Clearly $D_i^n \downarrow D_i$ as $n \rightarrow
\infty$ for $i=0,1,2$. Using the same arguments as in the proof of
Proposition 7 above we find that each set $D_i^n$ is convex for
$i=0,1,2$. Hence we can conclude that the boundary $b_i^n$ of
$D_i^n$ restricted to the square $[1/N,N]^2$ converges uniformly to
the boundary $b_i$ restricted to the square $[1/N,N]^2$ as $n
\rightarrow \infty$ for $i=0,1,2$. Thus, as in the case of the sets
$D_0, D_1, D_2$ and their boundaries $b_0, b_1, b_2$, the boundary
of the set $D_0^n$ restricted to the square $[1/N,N]^2$ is described
by a concave/continuous function $b_0^n : [1/n,\gamma_n] \rightarrow
[0,1]$ and the boundaries of the sets $D_1^n$ and $D_2^n$ are
described by a convex/continuous function $b_1^n : [1/N,N]
\rightarrow [1,N]$ for all $n \ge n_0$ with $n_0 \ge 1$ sufficiently
large.

\vspace{6pt}

3.\ We approximate the value function $\hat V$ by functions $\hat
V^n$ defined as follows
\begin{align} \h{5pc} \label{12.9}
\hat V^n(\varphi_1,\varphi_2) &= \hat V(\varphi_1,\varphi_2)\;\;
\text{if}\;\; (\varphi_1,\varphi_2) \in C_n \\[1pt] \notag &=
\hat M (\varphi_1,\varphi_2) \m \tfrac{1}{n}\;\; \text{if}\;\;
(\varphi_1,\varphi_2) \in D_n
\end{align}
for $(\varphi_1,\varphi_2) \in [0,\infty)^2$ with $n \ge n_0$ given
and fixed. Clearly $\hat V^n$ is a continuous function on
$[0,\infty)^2$ and moreover $\hat V^n$ restricted to $C_n$ and $D_n$
belongs to $C^2(\bar C_n)$ and $C^2(\bar D_n)$ respectively. Thus
the change-of-variable formula with local time on surfaces
\cite[Theorem 2.1]{Pe-2} is applicable to $\hat V^n$ composed with
$\varPhi^{\rho_N} = (\varPhi^{1,\rho_N},\varPhi^{2,\rho_N})$ and
this gives
\begin{align} \h{0pc} \label{12.10}
\hat V^n(\varPhi_t^{\rho_N}) &= \hat V^n(\varphi) + \int_0^t \hat
V_{\varphi_1}^n(\varPhi_s^{\rho_N})\, d \varPhi_s^{1,\rho_n} +
\int_0^t \hat V_{\varphi_2}^n(\varPhi_s^{\rho_N})\, d \varPhi_s
^{2,\rho_n} + \int_0^t \LL_{\varPhi^{\rho_N}} \hat V^n(\varPhi_s
^{\rho_N})\, ds \\ \notag &\h{-2pc}\h{13pt}+ \frac{1}{2} \int_0^t
\big[ \hat V_{\varphi_2}^n(\varPhi_s^{1,\rho_N},b_0^n(\varPhi_s^{1,
\rho_N})+) \m \hat V_{\varphi_2}^n(\varPhi_s^{1,\rho_N},b_0^n(\varPhi
_s^{1,\rho_N})-) \big]\, d \ell_s^{b_0^{1,n}}\!(\varPhi^{\rho_N})
\\ \notag &\h{-2pc}\h{13pt}+ \frac{1}{2} \int_0^t \big[\hat V_{
\varphi_1}^n(b_0^n(\varPhi_s^{2,\rho_N})+,\varPhi_s^{2,\rho_N}) \m
\hat V_{\varphi_1}^n(b_0^n(\varPhi_s^{2,\rho_N})-,\varPhi_s^{2,\rho_N})
\big]\, d \ell_s^{b_0^{2,n}}\!(\varPhi^{\rho_N}) \\ \notag &\h{-2pc}
\h{13pt}+ \frac{1}{2} \int_0^t \big[\hat V_{\varphi_2}^n(\varPhi_s
^{1,\rho_N},b_1^n(\varPhi_s^{1,\rho_N})+) \m \hat V_{\varphi_2}^n(
\varPhi_s^{1,\rho_N},b_1^n(\varPhi_s^{1,\rho_N})-) \big]\, d \ell
_s^{b_1^n}\!(\varPhi^{\rho_N}) \\ \notag &\h{-2pc}\h{13pt}+ \frac
{1}{2} \int_0^t \big[\hat V_{\varphi_1}^n(b_1^n(\varPhi_s^{2,\rho
_N})+,\varPhi_s^{2,\rho_N}) \m \hat V_{\varphi_1}^n(b_1^n(\varPhi
_s^{2,\rho_N})-,\varPhi_s^{2,\rho_N}) \big]\, d \ell_s^{b_2^n}\!
(\varPhi^{\rho_N}) \\ \notag &\h{-2pc}= \hat V^n(\varphi) + M_t
- \int_0^{t \wedge \rho_N}\!\! \hat H(\varPhi_s)\:\! I( \varPhi
_s\! \in\! C_n) \, ds + \frac{1}{2}\:\! L_t^{n,N}
\end{align}
for $\varphi \in [0,\infty)^2$ using \eqref{11.1} and \eqref{11.2}
where $M_t$ is a continuous martingale (the sum of the first two
integrals in the first identity of \eqref{12.10} above) and
$L_t^{n,N}$ is the sum of the final four integrals in the first
identity of \eqref{12.10} above. Note that the first partial
derivatives $\hat V_{\varphi_1}^n$ and $\hat V_{\varphi_2}^n$ are
discontinuous over the boundary curves $b_0^n$ and $b_1^n$ because
these boundary curves are not optimal. However, since $\hat V$ is
globally $C^1$ on $(0,\infty)^2$ by Corollary 15, it follows that
\begin{align} \h{4pc} \label{12.11}
&\sup_{1/N \le \varphi_1 \le N} \big \vert \hat V_{\varphi_2}^n(
\varphi_1, b_i^n(\varphi_1)) \m \hat M_{\varphi_2}^n(\varphi_1,b_i
^n (\varphi_1)) \big \vert \rightarrow 0 \\ \label{12.12} &\sup_{1/N
\le \varphi_1 \le N} \big \vert \hat V_{\varphi_1}^n(\varphi_2,b_i
^n(\varphi_2)) \m \hat M_{\varphi_1}^n(\varphi_2, b_i^n(\varphi_2))
\big \vert \rightarrow 0
\end{align}
for $i=0,1$ as $n \rightarrow \infty$. Note that the suprema in
\eqref{12.11} and \eqref{12.12} for $i=0,1$ provide uniform upper
bounds on the modulus of the integrands in the four integrals of
\eqref{12.10} with respect to the local times. To obtain a control
over the local times themselves in these four integrals (their
integrators) we now show that their expectations remain uniformly
bounded as the running time tends to infinity.

\vspace{6pt}

4.\ We first consider the case of $b_0^{1,n}$ and $b_0^{2,n}$
recalling that the two functions coincide by symmetry for $n \ge 1$
given and fixed. We thus focus on $b_0^{1,n}$ in the sequel. Since
$b_0^{1,n}$ is concave we see that $\varPhi^{2,\rho_N} -
b_0^{1,n}(\varPhi^{1,\rho_N})$ is a continuous semimartingale so
that by the It\^o-Tanaka formula we find that
\begin{align} \h{3pc} \label{12.13}
&\big( \varPhi_t^{2,\rho_N} \m b_0^{1,n}(\varPhi_t^{1,\rho_N})
\big)^+ = \big( \varphi_2 \m b_0^{1,n}(\varphi_1) \big)^+ \\
\notag &\h{14pt}+ \int_0^t I \big( \varPhi _s^{2,\rho_N}\! >\!
b_0^{1,n}(\varPhi_s^{1,\rho_N}) \big)\, d(\varPhi^{2,\rho_N} \m
b_0^{1,n}(\varPhi^{1,\rho_N}))_s + \frac{1}{2}\, \ell_t^{b_0^{1,
n}}\!\!(\varPhi^{\rho_N}) \\ \notag &= \big( \varphi_2 \m b_0^{1,
n}(\varphi_1) \big)^+ + \int_0^t I \big( \varPhi _s^{2,\rho_N}\!
>\! b_0^{1,n}(\varPhi_s^{1,\rho_N}) \big)\, d \varPhi_s^{2,\rho_N}
\\ \notag &\h{14pt}- \int_0^t I \big( \varPhi _s^{2,\rho_N}\! >\!
b_0^{1,n}(\varPhi_s^{1,\rho_N}) \big)\:\! (b_0^{1,n})'(\varPhi_s
^{1,\rho_N})\, d \varPhi_s^{1,\rho_N} \\ \notag &\h{14pt}- \int_
0^t I \big( \varPhi _s^{2,\rho_N}\! >\! b_0^{1,n}(\varPhi_s^{1,
\rho_N}) \big)\! \int_0^\infty d \ell_s^{\psi_1}\! (\varPhi^{1,
\rho_N})\, d (b_0^{1,n})'(\psi_1) + \frac{1}{2}\, \ell_t^{b_0
^{1,n}}\!\! (\varPhi^{\rho_N})
\end{align}
for $t \ge 0$ where $(b_0^{1,n})'$ denotes the first derivative of
$b_0^{1,n}$ (its existence follows by the implicit function theorem
since smooth fit fails at $b_0^{1,n}$ as pointed out above). Since
$b_0^{1,n}$ is concave we see that $d (b_0^{1,n})'$ defines a
non-positive measure on $[1/N,\gamma_n]$ so that the final integral
in \eqref{12.13} is non-positive. Using this fact in \eqref{12.13}
we obtain the following pathwise bound on the size of the local time
\begin{equation} \h{6pc} \label{12.14}
\ell_t^{b_0^{1,n}}\!\! (\varPhi^{\rho_N}) \le 2 \big( \varPhi_t^{2,
\rho_N} \m b_0^{1,n}(\varPhi_t^{1,\rho_N}) \big)^+\! - M_t
\end{equation}
where $M_t$ is a continuous martingale (the difference between the
second and the third integral in \eqref{12.13} above) for $t \ge 0$.
Taking $\EE_{\varphi_1,\varphi_2}^0$ on both sides of \eqref{12.14}
above and using that $\varPhi_t^{2,\rho_N} \le N$ for all $t \ge 0$
we find that
\begin{equation} \h{8pc} \label{12.15}
\EE_{\varphi_1,\varphi_2}^0 \big[ \ell_t^{b_0^{i,n}} \big] \le 2 N
\end{equation}
for all $t \ge 0$ and all $n \ge n_0$ with $(\varphi_1,\varphi_2)
\in [0,\infty)^2$ and $i=1,2$ (where the case $i=2$ follows from the
case $i=1$ by symmetry).

\vspace{6pt}

5.\ We next consider the case of $b_1^n$ and $b_2^n$ recalling that
the two functions coincide by symmetry for $n \ge 1$ given and
fixed. We thus focus on $b_1^n$ in the sequel. Similarly, since
$b_1^n$ is convex we see that $b_1^n(\varPhi^{1,\rho_N}) -
\varPhi^{2,\rho_N}$ is a continuous semimartingale so that by the
It\^o-Tanaka formula we find that
\begin{align} \h{3pc} \label{12.16}
&\big( b_1^n(\varPhi_t^{1,\rho_N}) \m \varPhi_t^{2,\rho_N}
\big)^+ = \big( b_1^n(\varphi_1) \m \varphi_2 \big)^+ \\
\notag &\h{14pt}+ \int_0^t I \big( b_1^n(\varPhi_s^{1,\rho_N})
\! >\! \varPhi _s^{2,\rho_N} \big)\, d(b_1^n(\varPhi^{1,\rho_N})
\m \varPhi^{2,\rho_N})_s + \frac{1}{2}\, \ell_t^{b_1^n}\! (\varPhi
^{\rho_N}) \\ \notag &= \big( b_1^n(\varphi_1) \m \varphi_2 \big)
^+ + \int_0^t I \big( b_1^n(\varPhi_s^{1,\rho_N})\! >\! \varPhi
_s^{2,\rho_N} \big)\:\! (b_1^n)'(\varPhi_s^{1,\rho_N})\, d \varPhi
_s^{1,\rho_N} \\ \notag &\h{14pt}+\int_0^t I \big( b_1^n(\varPhi
_s^{1,\rho_N})\! >\! \varPhi_s^{2,\rho_N} \big)\! \int_0^\infty
d\ell_s^{\psi_1}\! (\varPhi^{1,\rho_N})\, d (b_1^n)'(\psi_1)
\\ \notag &\h{14pt}- \int_0^t I \big( b_1^n(\varPhi_s^{1,
\rho_N})\! >\! \varPhi _s^{2,\rho_N} \big)\, d \varPhi_s^{2,
\rho_N} + \frac{1}{2}\, \ell_t^{b_1^n}\! (\varPhi^{\rho_N})
\end{align}
for $t \ge 0$ where $(b_1^n)'$ denotes the first derivative of
$b_1^n$ (its existence follows by the implicit function theorem
since smooth fit fails at $b_1^n$ as pointed out above). Since
$b_1^n$ is convex we see that $d (b_1^n)'$ defines a non-negative
measure on $[1/N,\gamma_n]$ so that the second last integral in
\eqref{12.16} is non-negative. Using this fact in \eqref{12.16} we
obtain the following pathwise bound on the size of the local time
\begin{equation} \h{6pc} \label{12.17}
\ell_t^{b_1^n}\! (\varPhi^{\rho_N}) \le 2 \big( b_1^n(\varPhi_t
^{1,\rho_N}) \m \varPhi_t^{2,\rho_N} \big)^+\!\! - M_t
\end{equation}
where $M_t$ is a continuous martingale (the difference between the
second and the final integral in \eqref{12.16} above) for $t \ge 0$.
Taking $\EE_{\varphi_1,\varphi_2}^0$ on both sides of \eqref{12.17}
above and using that $b_1^n \le b_1$ with $M_N := \sup_{\:\!1/N \le
\varphi_1 \le N} b_1(\varphi_1) < \infty$ we find that
\begin{equation} \h{8pc} \label{12.18}
\EE_{\varphi_1,\varphi_2}^0 \big[ \ell_t^{b_1^n} \big] \le 2 M_N
\end{equation}
for all $t \ge 0$ and all $n \ge n_0$ with $(\varphi_1,\varphi_2)
\in [0,\infty)^2$.

\vspace{6pt}

6.\ Combining \eqref{12.11}+\eqref{12.12} with
\eqref{12.15}+\eqref{12.18} we find that
$\EE_{\varphi_1,\varphi_2}^0 \big[ L_t^{n,N} \big] \rightarrow 0$ as
$n \rightarrow \infty$ for every $(\varphi_1,\varphi_2) \in
[0,\infty)^2$ and $N \ge 1$ given and fixed. Taking
$\EE_{\varphi_1,\varphi_2}^0$ on both sides of \eqref{12.10},
letting $n \rightarrow \infty$ and using the monotone convergence
theorem due to $\hat H \ge 0$, we obtain the following identity
\begin{equation} \h{3pc} \label{12.19}
\EE_{\varphi_1,\varphi_2}^0 \big[ \hat V( \varPhi_t^{\rho_N} ) \big]
\big] = \hat V(\varphi_1,\varphi_2) - \EE_{\varphi_1,\varphi_2}^0
\Big[ \int_0^{t \wedge \rho_N} \hat H(\varPhi_s)\:\! I(\varPhi_s
\! \in\! C)\, ds \Big]
\end{equation}
for $t \ge 0$ and $N \ge 1$ with $(\varphi_1,\varphi_2) \in
[0,\infty)^2$. Recalling that $0 \le \hat V \le \hat M$ where $\hat
M$ is defined in \eqref{4.13} above, and noting that
$\EE_{\varphi_1,\varphi_2}^0 \big( \sup_{0 \le s \le t} \varPhi_s^i
\big) < \infty$ for $i=1,2$, we see by letting $N \rightarrow
\infty$ that the dominated convergence theorem is applicable to the
left-hand side of \eqref{12.19}, while the monotone convergence
theorem is applicable to the right-hand side of \eqref{12.19} since
$\hat H \ge 0$. Letting $N \rightarrow \infty$ in \eqref{12.19} we
thus obtain
\begin{equation} \h{3pc} \label{12.20}
\EE_{\varphi_1,\varphi_2}^0 \big[ \hat V( \varPhi_t ) \big]
\big] = \hat V(\varphi_1,\varphi_2) - \EE_{\varphi_1,\varphi_2}^0
\Big[ \int_0^t \hat H(\varPhi_s)\:\! I(\varPhi_s\! \in\! C)\, ds \Big]
\end{equation}
for $t \ge 0$ and $(\varphi_1,\varphi_2) \in [0,\infty)^2$.

\vspace{6pt}

7.\ Despite the fact that neither $(\varPhi_t^1)_{t \ge 0}$ nor
$(\varPhi_t^2)_{t \ge 0}$ is uniformly integrable (since
$\varPhi_t^i \rightarrow 0$ with
$\PP_{\!\varphi_1,\varphi_2}$\!-probability one as $t \rightarrow
\infty$ but $\EE_{\varphi_1,\varphi_2}(\varPhi_t^i) = \varphi_i$ for
all $t \ge 0$ with $i=1,2$ and $(\varphi_1,\varphi_2) \in
(0,\infty)^2$ given and fixed) we claim that
\begin{equation} \h{5pc} \label{12.21}
\{\:\! \hat M(\varPhi_t^1,\varPhi_t^2)\; \vert\; t \ge 0\:\! \}\;\;
\text{is uniformly integrable}
\end{equation}
where we recall that $\hat M$ is defined in \eqref{4.13} above. For
this, note that $0 \le \hat M(\varPhi_t^1,\varPhi_t^2) = c\;\! \big(
(\varPhi_t^1 \p \varPhi_t^2) \wedge (1 \p \varPhi_t^1) \wedge (1 \p
\varPhi_t^2) \big) \le c\;\! \big( (1 \p \varPhi_t^1) \wedge (1 \p
\varPhi_t^2) \big) = c\;\! \big( 1 \p \varPhi_t^1 \wedge \varPhi_t^2
\big)$ for $t \ge 0$. A direct martingale argument based on
\eqref{5.8} then gives
\begin{align} \h{3pc} \label{12.22}
\EE_{\varphi_1,\varphi_2}^0 \big[ \varPhi_t^1 \wedge \varPhi_t^2 \big]
&= \varphi_1\:\! \varphi_2\, \EE \big[ e^{\frac{\mu}{\sqrt{2}} \sqrt{3}
\;\! W_t^1 - \mu^2t} \big( e^{\frac{\mu}{\sqrt{2}}\;\! W_t^2} \wedge
e^{-\frac{\mu}{\sqrt{2}}\;\! W_t^2} \big) \big] \\[3pt] \notag &\le
\varphi_1\:\! \varphi_2\, e^{-\frac{1}{4} \mu^2 t}\, \EE \big[ e^{
\frac{\mu}{\sqrt{2}} \sqrt{3} \;\! W_t^1 - \frac{3}{4} \mu^2t} \big]
= \varphi_1\:\! \varphi_2\, e^{-\frac{1}{4} \mu^2 t} \rightarrow 0
\end{align}
as $t \rightarrow \infty$ for $(\varphi_1,\varphi_2) \in
[0,\infty)^2$. Since $\varPhi_t^1 \wedge \varPhi_t^2 \rightarrow 0$
with $\PP_{\!\varphi_1,\varphi_2}$\!-probability one as $t
\rightarrow \infty$ for $(\varphi_1,\varphi_2) \in (0,\infty)^2$
given and fixed, we see from \eqref{12.22} that $\{\;\! \varPhi_t^1
\wedge \varPhi_t^2\; \vert\; t \ge 0\;\! \}$ is uniformly
integrable. Hence by the bound preceding to \eqref{12.22} we see
that \eqref{12.21} holds as claimed.

\vspace{6pt}

8.\ Since $0 \le \hat V(\varPhi_t^1,\varPhi_t^2) \le \hat
M(\varPhi_t^1,\varPhi_t^2)$ for $t \ge 0$ we see from \eqref{12.21}
that $\{\;\! \hat V(\varPhi_t^1,\varPhi_t^2)\; \vert\; t \ge 0\;\!
\}$ is uniformly integrable. Letting $t \rightarrow \infty$ in
\eqref{12.20} and using that $\hat V(\varPhi_t^1,\varPhi_t^2)
\rightarrow 0$ with $\PP_{\!\varphi_1,\varphi_2}$\!-probability one
we thus find by the extended dominated convergence theorem (applied
to the left-hand side) and the monotone convergence theorem (applied
to the right-hand side) that the following identity holds
\begin{equation} \h{5pc} \label{12.23}
\hat V(\varphi_1,\varphi_2) = \EE_{\varphi_1,\varphi_2}^0 \Big[ \int
_0^\infty\! \hat H(\varPhi_s)\:\! I(\varPhi_s\! \in\! C)\, ds \Big]
\end{equation}
$(\varphi_1,\varphi_2) \in [0,\infty)^2$. Combining \eqref{12.23}
with \eqref{12.3} we obtain \eqref{12.6} as claimed. Evaluating
$\hat V$ from \eqref{12.23} at the optimal stopping points
$(\varphi_1,b_0(\varphi_1))$ and $(\varphi_1,b_1(\varphi_1))$ upon
using that $\hat V(\varphi_1,b_0(\varphi_1)) = \hat
M(\varphi_1,b_0(\varphi_1)) = \varphi_1 \p b_0(\varphi_1)$ and $\hat
V(\varphi_1,b_1(\varphi_1)) = \hat M(\varphi_1,b_1(\varphi_1)) = 1
\p \varphi_1$ for $\varphi_1 \in [0,\gamma]$ and $\varphi_1 \in
[0,\infty)$ respectively, we see that the functions $b_0$ and $b_1$
solve the integral equations \eqref{12.4} and \eqref{12.5} as
claimed. This completes the proof of the existence of the solution
to these equations.

\vspace{6pt}

(II) \emph{Uniqueness}. To show that $b_0$ and $b_1$ are a unique
solution to the system \eqref{12.4}+\eqref{12.5} one can adopt the
four-step procedure from the proof of uniqueness given in
\cite[Theorem 4.1]{DP} extending and further refining the original
arguments from \cite[Theorem 3.1]{Pe-1} in the case of a single
boundary. Given that the present setting creates no additional
difficulties we will omit further details of this verification and
this completes the proof. \hfill $\square$

\vspace{12pt}

The coupled system of nonlinear Fredholm integral equations
\eqref{12.4}+\eqref{12.5} can be used to find the optimal stopping
boundaries $b_0$ and $b_1$ numerically (using Picard iteration).
Inserting these $b_0$ and $b_1$ into \eqref{12.6} we also obtain a
closed form expression for the value function $\hat V$. Collecting
the results derived throughout we now disclose the solution to the
initial problem.

\vspace{12pt}

\textbf{Corollary 20.} \emph{The value function of the initial problem
\eqref{3.3} is given by
\begin{equation} \h{7pc} \label{12.24}
V(\pi_0,\pi_1,\pi_2) = \pi_0\;\! \hat V \Big( \frac{\pi_1}{\pi_0},
\frac{\pi_2} {\pi_0} \Big)
\end{equation}
for $(\pi_0,\pi_1,\pi_2) \in [0,1]^3$ with $\pi_0 \p \pi_1 \p \pi_2
= 1$ where the function $\hat V$ is given by \eqref{12.6} above. The
optimal stopping time in the initial problem \eqref{3.3} is given by
\begin{align} \h{0pc} \label{12.25}
\tau_* = \inf\, \{\, t \ge 0\; &\vert\; \tfrac{\pi_i}{\pi_0}\:\!
e^{\mu (X_t^i - X_t^0)} \ge \tfrac{\pi_j}{\pi_0}\:\! e^{\mu (X_t^j -
X_t^0)} \;\; \text{with}\;\; \tfrac{\pi_i}{\pi_0}\:\! e^{\mu (X_t^i
- X_t^0)} \le b_0(\tfrac{\pi_j}{\pi_0}\:\! e^{\mu (X_t^j - X_t^0)})
\;\; \\[2pt] \notag &\h{8pt}\text{or}\;\; \tfrac{\pi_i}{\pi_0}\:\!
e^{\mu (X_t^i - X_t^0)} \ge b_1(\tfrac{\pi_j}{\pi_0}\:\! e^{\mu
(X_t^j - X_t^0)}) \;\; \text{for}\;\; i \ne j\;\; \text{in} \;\;
\{1,2\}\, \}
\end{align}
where $b_0$ and $b_1$ are a unique solution to the coupled system of
nonlinear Fredholm integral equations \eqref{12.4}+\eqref{12.5}. The
optimal decision function $d_{\tau_*}$ in the initial problem
\eqref{3.3} equals $0$ if stopping in \eqref{12.25} happens at
$b_0$, equals $1$ if stopping in \eqref{12.25} happens at $b_1$ with
$i=1$, and equals $2$ is stopping in \eqref{12.25} happens at $b_1$
with $i=2$.}

\vspace{12pt}

\textbf{Proof.} The identity \eqref{12.24} was established in
\eqref{4.11} above. The explicit form \eqref{12.25} follows from
\eqref{12.7} in Theorem 19 combined with \eqref{4.2}-\eqref{4.4}
above. The final claim on the optimal decision function follows from
\eqref{3.7} combined with the argument used in the second equality
of \eqref{4.15} above completing the proof. \hfill $\square$

\vspace{6pt}

\begin{center}

\end{center}

%%%%%%%%%%%%%%%%%%%%%%%%%%%%%%%%%%%%%%%%%%%%%%%%%%%%%%%%%%%%%%%%%%%%%%%%%%%%%%%
%%% Affiliations %%%
%%%%%%%%%%%%%%%%%%%%%%%%%%%%%%%%%%%%%%%%%%%%%%%%%%%%%%%%%%%%%%%%%%%%%%%%%%%%%%%

\par \leftskip=24pt

\vspace{20pt}

\ni Philip A.\ Ernst \\
Department of Statistics \\
Rice University \\
6100 Main Street \\
Houston TX 77005 \\
United States \\
\texttt{philip.ernst@rice.edu}

\leftskip=15pc \vspace{-101pt}

\ni Goran Peskir \\
School of Mathematics \\
The University of Manchester \\
Oxford Road \\
Manchester M13 9PL \\
United Kingdom \\
\texttt{goran@maths.man.ac.uk}

\leftskip=29.5pc \vspace{-102pt}

\ni Quan Zhou \\
Department of Statistics \\
Rice University \\
6100 Main Street \\
Houston TX 77005 \\
United States \\
\texttt{qz9@rice.edu}

\par

\end{document}